\documentclass[10pt]{nm}
\usepackage{graphicx}
\usepackage{subfigure}
\usepackage{algorithm}
\usepackage{algorithmic}
\newcommand{\tabincell}[2]{\begin{tabular}{@{}#1@{}}#2\end{tabular}}

\begin{document}

\markboth{L. Wang, Y. Huang and K. Jiang}{Error analysis of SAV FEM to PFC model}
\title{Error analysis of SAV finite element method to phase field crystal model}

\author[L. Wang, Y. Huang and K. Jiang]{Liupeng Wang\affil{1}, Yunqing Huang\affil{1}\comma\affil{2}\comma\affil{3} and  Kai Jiang\affil{1}\comma\affil{2}\comma\affil{3}\comma\corrauth}

\address{\affilnum{1}\ School of Mathematics and Computational Science, Xiangtan University, Xiangtan 411105, Hunan, P.R.China. \\
          \affilnum{2}\ Hunan Key Laboratory for Computation and Simulation in Science and Engineering, Xiangtan 411105, Hunan, P.R.China.\\
	\affilnum{3}\ Key Laboratory of Ministry of Education for Intelligent Computing and Information Processing, Xiangtan University, Xiangtan 411105, Hunan, P.R.China}

\emails{{\tt kaijiang@xtu.edu.cn} (K. Jiang)}

\begin{abstract}
In this paper, we construct and analyze an energy stable scheme
by combining the latest developed scalar auxiliary variable (SAV) approach
and linear finite element method (FEM) for phase field crystal
(PFC) model, and show rigorously that the scheme is first-order
in time and second-order in space for the $L^2$ and $H^{-1}$ gradient flow equations. To
reduce efficiently computational cost and capture accurately the
phase interface, we give a simple adaptive strategy, equipped
with a posteriori gradient estimator, \textit{i.e.} $L^2$ norm of
the recovered gradient. Extensive numerical experiments are
presented to verify our theoretical results and to demonstrate
the effectiveness and accuracy of our proposed method.
\end{abstract}



\maketitle

\section{Introduction}\label{secInt}
The phase field crystal (PFC) model\cite{elder2002modeling,
elder2004modeling} was proposed as an approach
to simulate crystals at the atomic scale but on a
coarse-grained diffusive time scale\cite{wise2009an}. 
Many physical processes, such as the formation of ordered structures,
phase separation of polynary systems, can be described using this
model. The PFC model can also explain elastic and plastic
deformations of the lattice, dislocations, grain boundaries,
multiple crystal orientations and many other observable
phenomena\cite{wise2009an, provatas2007using}. 

There are several kinds of PFC models. In general, they can be
classified into two classes according to characteristic length
scale: one-length-scale and multi-length-scale.
One-length-scale PFC models can be used to describe the phase
behavior of periodic structures\cite{brazovskii1975phase, swift1977hydrodynamic, ohta1986equilibrium}. Accordingly,
multi-length-scale PFC models can be employed to explain the
formation of quasicrystals\cite{lifshitz1997theoretical, jiang2014numerical}. 
In this work, we focus on
the development of numerical methods of one-length-scale PFC model.
In particular, the classic Landau-Brazovskii (LB)
model\cite{brazovskii1975phase, fredrickson1987fluctuation, zhang2008an}
will be used to demonstrate our proposed method.
The LB model was built to investigate the
character of phase transition. It has been discovered in many
scientific fields. For example, the LB model can be derived
from more complicated self-consistent field theory of diblock
copolymers\cite{fredrickson2006the}. Compared with the typical
Swift-Hohenberg (SH) model with double-well bulk
energy\,\cite{swift1977hydrodynamic}, LB energy functional
includes a cubic term which can be used to study the first-order
phase transition.

The $L^2$ (Allen-Cahn) or $H^{-1}$ (Cahn-Hilliard) gradient flow
equation is usually adopted to describe the dynamic behavior of the
phase-field or PFC model. These dynamic equations are
time-dependent nonlinear partial differential equations (PDEs). 
It is hard to find non-trivial analytical solutions. Therefore, 
numerically solving these nonlinear PDEs is an efficient approach.
To guarantee convergence, numerical schemes of these equations 
are required to satisfy the energy dissipation property. 
Meanwhile, an accurate and efficient approach should be designed
to deal with nonlinear terms. In terms of time discretization, 
there have been several
effective methods which can preserve energy dissipation law,
including the convex splitting method\cite{wise2009an, eyre1998an, 
baskaran2013convergence, vignal2015energy,li2017second,guo2019efficient,lee2019energy,
lee2019effective}, 
stabilized approach\cite{shen2010numerical,tang2016implicit,li2017on,yan2018second}, 
invariant energy quadratization (IEQ) method\cite{zhao2017numerical,guo2019efficient}
and recently developed scalar auxiliary variable (SAV) approach\cite{shen2018scalar}. By introduing a scale auxiliary variable to the nonlinear part of energy functional, the SAV approach has a modified energy dissipation property for a large class of gradient flows. 
The convergent and error analysis of semi-discrete SAV scheme
has been given by Shen and Xu\cite{shen2018convergence}.
The analysis of energy stability and convergence of fully
discretized SAV block-centered finite difference method has been
established for gradient flows\cite{li2019energy}. 
More studies about the PFC problem can be found in recent
literature\,\cite{baskaran2013convergence,grasselli2015energy,vignal2015energy,guo2016local,guo2018high}.

In the study to the PFC model, finite difference methods\cite{elder2004modeling,elder2007phase,wise2009an,hu2009stable,church2019high} or spectral methods\cite{cheng2006an,tegze2009advanced,shen2018scalar} are limited to regular regions, such as two-dimensional square region or three-dimensional cube region. For complex geometries, finite element method (FEM)\cite{backofen2007nucleation,praetorius2015efficient,diegel2019a} is a better choice. Furthermore, the FEM can be further combined with adaptive technologies, which are well suitable for the phase behavior of PFC models, such as the formation of ordered structures, phase transition processes, and coarse-grained processes. The adaptive method can effectively decrease the cost of computing and accurately capture the phase interface.

In this work, we will combine SAV time discretization and
FEM spatial discretization to solve the gradient flow
equation of LB model. Based on the energy dissipation and the
SAV scheme, the derivation process of $H^2$ bounds of the
solution is shown in detail. For our fully discrete scheme, 
we demonstrate its energy stability, and carry out error estimate.
Applying our method, we can effectively simulate the mesoscale
self-assembly of the diblock copolymer system in two-dimensional convex geometries. 
In addition, we will consider an adaptive FEM for the PFC model.
There are many adaptive finite element methods for phase field
equation\cite{feng2005posteriori,du2008adaptive,hu2009multi,wu2017space,chen2019scr}.
To reduce computational cost, we first apply an adaptive method which is
effective against phase field equation to the PFC model. 
Numerical results demonstrate that directly using the gradient as
the indicator is more efficient than the posterior error
estimator does in solving this problem. Since the gradient
obtained from the numerical solution may be discontinuous, 
therefore, a smooth recovered gradient is employed as the adaptive
indicator in our adaptive FEM.   

The rest of this paper is organized as follows. In
Sec.\,\ref{sec:DE}, we introduce the LB free energy functional
and take the $L^2$ gradient flow as an example to 
derive its Allen-Cahn dynamic equation. Sec.\,\ref{sec:NM}
details our numerical method, which consists of discretization
schemes, energy dissipation and error estimate for $L^2$ and
$H^{-1}$ dynamical equations. In
Sec.\,\ref{sec:NR}, numerical experiments are given to illustrate
the accuracy and effectiveness of our scheme.  
Several standard ordered structures on two-dimensional convex regions
are also obtained in this section. 
Sec.\,\ref{sec:AFEM} gives a simple but efficient adaptive FEM to
PFC model. In Sec.\,\ref{sec:CO}, we give conclusions and
outlooks. 

\section{Physical model}
\label{sec:DE}

The dimensionless free energy functional of LB model is\cite{zhang2008an} 
\begin{equation}
\mathcal{E}(\phi(\mathbf{r})) = \int_{\Omega} \left\{
\frac{\xi^2}{2} [(\nabla^2+1)\phi(\mathbf{r})]^2 +
\frac{\alpha}{2}[\phi(\mathbf{r})]^2 -
\frac{\gamma}{3!}[\phi(\mathbf{r})]^3 +
\frac{1}{4!}[\phi(\mathbf{r})]^4 \right\} \mathrm d\mathbf{r}.
\end{equation}
where $\Omega$ is a two dimensional bounded domain with
Lipschitz boundary $\partial \Omega$, $\phi$ is the density
deviation of a kind of monomer from the disordered phase, 
$\xi, \alpha$ and $\gamma$
are the parameters of the model, $\nabla^2$ is the Laplace operator. 

The gradient flow of LB model is
\begin{equation}
	\phi_t = \mathcal{L}\frac{\delta \mathcal{E}}{\delta \phi}, \label{eq:MAC}
\end{equation}
where $\mathcal{L}$ is a negative operator. For $L^2$ gradient
flow $\mathcal{L}=-\mathcal{I}$, while for $H^{-1}$ gradient flow,
$\mathcal{L}=\nabla^2$.
For brevity, in the following, we take the $L^2$ gradient flow
as an example to derive the Allen-Cahn dynamical equation. 
The $H^{-1}$ gradient flow equation can be derived similarly.

Now we will give the boundary conditions to guarantee the energy
dissipation property of the Allen-Cahn dynamical equation. Denote that
\begin{equation}\label{eq:nonlinearderv}
	\mathcal{N}(\phi):=\frac{\alpha}{2}\phi^2 - \frac{\gamma}{3!}\phi^3 + \frac{1}{4!}\phi^4,
\end{equation}
then
\begin{align}
\mathcal{N}^{\prime}(\phi)=\alpha\phi - \frac{\gamma}{2}\phi^2
	+ \frac{1}{3!}\phi^3, ~~~
\mathcal{N}''(\phi)=\alpha - \gamma\phi + \frac{1}{2}\phi^2, 
\end{align}
therefore, the LB model becomes 
\begin{equation}
\mathcal{E}(\phi(\mathbf{r})) = \int_{\Omega} \left\{
\frac{\xi^2}{2} [(\nabla^2+1)\phi]^2 +\mathcal{N}(\phi) \right\} \mathrm d\mathbf{r}.
\end{equation}

The free energy $\mathcal{E}$ take a derivative with respect to time $t$ is 
\begin{eqnarray*}
\frac{\partial \mathcal{E}(\phi)}{\partial t} &=& \int_{\Omega} \{ \xi^2 [(\nabla^2+1)\phi(\nabla^2+1) \phi_t]  + \mathcal{N}^{\prime}(\phi) \phi_t \} \mathrm d\mathbf{r} \\
&=& \int_{\Omega}  [\xi^2 (\nabla^2+1)\phi +\mathcal{N}^{\prime}(\phi)] \phi_t \mathrm d\mathbf{r}  +  \xi^2 \int_{\Omega}(\nabla^2+1)\phi  \nabla^2 \phi_t  \mathrm d\mathbf{r} \\
&=&   \int_{\Omega}  [\xi^2 (\nabla^2+1)\phi +\mathcal{N}^{\prime}(\phi)] \phi_t \mathrm d\mathbf{r} + \xi^2 \int_{\partial \Omega}  (\nabla^2+1)\phi\, \nabla \phi_t \cdot \mathbf{n} \mathrm d S  \\
&&- \xi^2 \int_{\Omega}\nabla (\nabla^2+1)\phi \cdot \nabla \phi_t  \mathrm d\mathbf{r} \\
&=& \int_{\Omega}  [\xi^2 (\nabla^2+1)\phi +\mathcal{N}^{\prime}(\phi)] \phi_t \mathrm d\mathbf{r} + \xi^2 \int_{\partial \Omega}  (\nabla^2+1)\phi\, \nabla \phi_t \cdot \mathbf{n} \mathrm d S \\
&& - \xi^2 \int_{\partial \Omega} \nabla (\nabla^2+1)\phi \cdot \mathbf{n}\, \phi_t  \mathrm d S+ \xi^2\int_{\Omega}\nabla^2 (\nabla^2+1)\phi \, \phi_t  \mathrm d\mathbf{r},  
\end{eqnarray*}
here we introduce two Neumann boundary conditions 
\begin{equation}\label{eq:NBC}
\nabla \phi \cdot \mathbf{n}|_{\partial \Omega}=0, \quad  \nabla (\nabla^2+1)\phi \cdot \mathbf{n} |_{\partial \Omega}=0,
\end{equation}
then
\begin{equation}
\frac{\partial \mathcal{E}(\phi)}{\partial t}= \int_{\Omega}
[\xi^2 (\nabla^2+1)^2\phi + \mathcal{N}'(\phi)] \phi_t \mathrm d\mathbf{r}.
\label{eq:BT}
\end{equation}

Define a function space $W(\Omega)$ as
\begin{equation}
W(\Omega) := \{w \in H^2(\Omega): \nabla w \cdot \mathbf{n}|_{\partial \Omega}=0
,\nabla (\nabla^2+1)w \cdot \mathbf{n} |_{\partial \Omega}=0\}.
\end{equation}
In the sense of the Gateaux differential for all $v \in W(\Omega)$, we have the following equation
\begin{eqnarray*}
\left( \frac{\delta \mathcal{E}}{\delta \phi}, v
\right) &=& \frac{d}{d \theta} \mathcal{E}(\phi+\theta v) |_{\theta=0}\\
&=&  \lim\limits_{\theta \to 0}\frac{1}{\theta}[\mathcal{E}(\phi+\theta v) - \mathcal{E}(\phi)]\\
&=& \big(\xi^2 (\nabla^2+1)^2\phi +\mathcal{N}'(\phi),v\big), 
\end{eqnarray*}
where $\phi \in W(\Omega)$, $(\cdot,\cdot)$ denote the $L^2$ inner product.
According to the variational principle, we have
\begin{equation}
\frac{\delta \mathcal{E}}{\delta \phi} = \xi^2 (\nabla^2+1)^2\phi
+\mathcal{N}'(\phi).
\label{eq:Bphi}
\end{equation}
From Eqns.\,\eqref{eq:MAC}, \eqref{eq:BT} and \eqref{eq:Bphi}, it is easy to verify
\begin{equation}
\frac{\partial \mathcal{E}(\phi)}{\partial t}= - \int_{\Omega}  [\xi^2 (\nabla^2+1)^2\phi +\mathcal{N}'(\phi)]^2 \mathrm d\mathbf{r} \le 0.
\label{eq:BTT}
\end{equation}
Therefore, Allen-Cahn equation \eqref{eq:MAC} satisfies energy
dissipation with the Neumann boundary conditons \eqref{eq:NBC}.

Combining Eqns.\,\eqref{eq:MAC}, \eqref{eq:nonlinearderv},
\eqref{eq:NBC} and \eqref{eq:Bphi}, 
the governing equation can be written as
\begin{subequations}\label{eq:LB-AC}
\begin{equation}
\phi_t = -\xi^2 (\nabla^2+1)^2\phi-\mathcal{N}'(\phi),
\end{equation}
\begin{equation}
\nabla \phi \cdot \mathbf{n}|_{\partial \Omega}=0, \quad \nabla (\nabla^2+1)\phi \cdot \mathbf{n} |_{\partial \Omega}=0.
\end{equation}
\end{subequations}

By introducing a new function $\psi = (\nabla^2+1) \phi$, we
can write Eqns.\,\eqref{eq:LB-AC} as the coupled system

\begin{subequations}\label{eq:LB-AC-SP}
\begin{equation}
\phi_t = -\xi^2(\nabla^2+1) \psi - \mathcal{N}'(\phi),
\end{equation}
\begin{equation}
\psi = (\nabla^2+1) \phi,
\end{equation}
\begin{equation}
\nabla \phi \cdot \mathbf{n} |_{\partial \Omega}=0, \quad  \nabla \psi \cdot \mathbf{n} |_{\partial \Omega}=0.
\end{equation}
\end{subequations} 
 
\textbf{Remark 2.1}\quad  The splitting technique used in Eqns.\,\eqref{eq:LB-AC-SP} is valid
for convex regions\cite{brenner2011c0}.

\section{Numerical Methods}
\label{sec:NM}

The main aim of this section is to propose numerical methods
to solve the gradient flows of LB model.
For brevity, for the Allen-Cahn equation, we first present the
discretization scheme, prove the energy stability, and give the
error estimate in detail. Then the corresponding results about
the Cahn-Hilliard equation are also given.

\subsection{SAV discretization}
Recently, Shen \textit{et al.}
proposed an efficient time discretization scheme,
\textit{i.e.}, the SAV scheme, to a class of gradient flows\cite{shen2018scalar}. 
Here, we shall apply the idea of the SAV approach to discretize
Eqns.\,\eqref{eq:LB-AC-SP} in time direction.

Let $\mathcal{E}_1(\phi)=\int_{\Omega} \mathcal{N}(\phi) \mathrm
d\mathbf{r} $. Then we introduce the scalar auxiliary variable
$s=\sqrt{\mathcal{E}_1(\phi)+D_0}$, where $D_0$ is a constant to ensure $\mathcal{E}_1(\phi)+D_0 \ge 0$, and write Eqns.\,\eqref{eq:LB-AC-SP} as 
\begin{subequations}\label{eq:LB-AC-SAV}
\begin{equation}
\phi_t = -\left\{\xi^2(\nabla^2+1) \psi + u(\phi)s\right\} ,\label{eq:LB-AC-SAV1}
\end{equation}
\begin{equation}
\psi = (\nabla^2+1) \phi,\label{eq:LB-AC-SAV2}
\end{equation}
\begin{equation}
s_t = \frac{1}{2}\int_{\Omega} u(\phi) \phi_t \mathrm d\mathbf{r},\label{eq:LB-AC-SAV3}
\end{equation}
\end{subequations} 
 where $u(\phi):=\frac{\mathcal{N}'(\phi)}{\sqrt{\mathcal{E}_1(\phi)+D_0}}$.

Assume $\Delta t$ is a fixed time step and $\phi^{n}$ is the approximation of $\phi(\mathbf{r},t^{n})$ at time $t^{n}=n \Delta t$, we can construct the first-order SAV scheme:
\begin{subequations}\label{eq:LB-AC-RFOSAV}
\begin{align}
\frac{\phi^{n+1}-\phi^{n}}{\Delta t} &= -\{\xi^2(\nabla^2+1) \psi^{n+1} +u(\phi^{n}) s^{n+1} \},\\
\psi^{n+1} &= (\nabla^2+1) \phi^{n+1},\\
s^{n+1}-s^{n} &=  \frac{1}{2} \big(u(\phi^{n}), \phi^{n+1}-\phi^{n}\big).
\end{align}
\end{subequations}

\textbf{Remark 3.1}\quad  The higher-order schemes based on SAV technique
	can be easily constructed, see Ref.\cite{shen2018scalar} for more details.

\subsection{FEM discretization}
We discretize Eqns.\,\eqref{eq:LB-AC-RFOSAV} in space using the FEM.
Let $V(\Omega)$ denote both the trial and test function spaces
\begin{equation}
V(\Omega):= \{w \in H^1(\Omega):\nabla w \cdot \mathbf{n}|_{\partial \Omega}=0\}.
\end{equation} 
The corresponding Galerkin form of Eqns.\,\eqref{eq:LB-AC-RFOSAV} can be stated as follows:
for $\forall v \in V$, find $\phi,\psi \in V$ such that:
\begin{subequations}\label{eq:LB-AC-RFO}
\begin{align}
\left(\frac{\phi^{n+1}-\phi^{n}}{\Delta t},v\right) &= -\{\xi^2\big(\psi^{n+1},v\big)-\xi^2\big(\nabla \psi^{n+1},\nabla v\big) +\big(u(\phi^{n}),v\big) s^{n+1} \},\\
\big(\psi^{n+1},v\big) &= \big(\phi^{n+1},v\big)-\big(\nabla \phi^{n+1},\nabla v\big),\\
s^{n+1}-s^{n} &= \frac{1}{2}\big(u(\phi^{n}),\phi^{n+1}-\phi^{n}\big).
\end{align}
\end{subequations}

Let $\mathcal{T}_h$ be a comforming mesh of $\Omega$ with $h=\max\limits_{k\in \mathcal{T}_h}\{h_{\tau}\}$, $\tau$ be line segment in 1D or triangle in 2D, and $V_h$ be the linear finite element space over $\mathcal{T}_h$ defined by
$$V_h :=\{w \in H^1(\Omega): w|_{\tau} \in P_1(\tau), \forall k \in \mathcal{T}_h\}.$$
Thus, Eqns.\,\eqref{eq:LB-AC-RFO} is transformed as follow: find $\phi_h,\psi_h \in V_h$, such that for $\forall  v_h \in V_h$:
\begin{subequations}\label{eq:LB-AC-RFOSAV-LFEM}
\begin{align}
\left(\frac{\phi_h^{n+1}-\phi_h^{n}}{\Delta t},v_h\right) &= -\{\xi^2\big(\psi_h^{n+1},v_h\big)-\xi^2\big(\nabla \psi_h^{n+1},\nabla v_h\big)  +\big(u(\phi_h^n),v_h\big) s_h^{n+1}\},\label{eq:LB-AC-RFOSAV-LFEM1}\\
\big(\psi_h^{n+1},v_h\big) &= \big(\phi_h^{n+1},v_h\big)-\big(\nabla \phi_h^{n+1},\nabla v_h\big),\label{eq:LB-AC-RFOSAV-LFEM2}\\
s_h^{n+1}-s_h^{n} &=  \frac{1}{2}\big(u(\phi_h^n),\phi_h^{n+1}-\phi_h^{n}\big)\label{eq:LB-AC-RFOSAV-LFEM3},
\end{align}
\end{subequations}

\subsection{Energy stability}

We shall prove the energy stability of Eqns.\,\eqref{eq:LB-AC-RFOSAV-LFEM}. The norm of $L^2(\Omega)$ is denoted by $\|\cdot\|$.

\textbf{Theorem 3.1}\quad  If we denote the modified energy 
$$\tilde{\mathcal{E}}(\psi,s) := \frac{\xi^2}{2}\|\psi\|^2+s^2 ,$$
then Eqns.\,\eqref{eq:LB-AC-RFOSAV-LFEM} is unconditionally energy stable with the modified energy.

\textbf{Proof:}\quad We take $v_h=\phi_h^{n+1}-\phi_h^{n}$ in Eqn.\,\eqref{eq:LB-AC-RFOSAV-LFEM1} , and find
\begin{equation}\label{eq:pr1}
	\begin{aligned}
\frac{1}{\Delta t}\|\phi_h^{n+1}-\phi_h^{n}\|^2+\xi^2\big(\psi_h^{n+1},\phi_h^{n+1}-\phi_h^{n}\big)-\xi^2\big(\nabla \psi_h^{n+1},\nabla (\phi_h^{n+1}-\phi_h^{n})\big)\\ +\big(u(\phi_h^{n}),\phi_h^{n+1}-\phi_h^{n}\big) s_h^{n+1}=0 .
\end{aligned}
\end{equation}
According to Eqn.\,\eqref{eq:LB-AC-RFOSAV-LFEM2}, we have
\begin{equation}\label{eq:ns1}
\big(\psi_h^{n},v_h\big) = \big(\phi_h^{n},v_h\big)-\big(\nabla \phi_h^{n},\nabla v_h\big).
\end{equation}
Substracting Eqn.\,\eqref{eq:LB-AC-RFOSAV-LFEM2} by Eqn.\,\eqref{eq:ns1}, we have
\begin{equation}\label{eq:ns2}
\big(\psi_h^{n+1}-\psi_h^{n},v_h\big) = \big(\phi_h^{n+1}-\phi_h^{n},v_h\big)-\big(\nabla (\phi_h^{n+1}-\phi_h^{n}),\nabla v_h\big).
\end{equation}
Setting $v_h=\psi_h^{n+1}$ to Eqn.\,\eqref{eq:ns2} and using the identity
$$(a-b,2a)=|a|^2-|b|^2+|a-b|^2,$$ we have
\begin{equation}\label{eq:pr2}
\begin{aligned}
&\big(\phi_h^{n+1}-\phi_h^{n},\psi_h^{n+1}\big)-\big(\nabla (\phi_h^{n+1}-\phi_h^{n}),\nabla \psi_h^{n+1}\big){}\\
&= \big(\psi_h^{n+1}-\psi_h^{n},\psi_h^{n+1}\big){}\\
&= \frac{1}{2}(\|\psi_h^{n+1}\|^2-\|\psi_h^{n}\|^2+\|\psi_h^{n+1}-\psi_h^{n}\|^2).
\end{aligned}
\end{equation}
Multiplying Eqn.\,\eqref{eq:LB-AC-RFOSAV-LFEM3} with $2s_h^{n+1}$, we obtain
\begin{equation}\label{eq:pr3}
\begin{aligned}
&\big(u(\phi_h^{n}),\phi_h^{n+1}-\phi_h^{n}\big)s_h^{n+1} \\
&= 2(s_h^{n+1}-s_h^{n})s_h^{n+1}{}\\
&=  (s_h^{n+1})^2-(s_h^{n})^2+(s_h^{n+1}-s_h^{n})^2.
\end{aligned}
\end{equation}
Substituting Eqn.\,\eqref{eq:pr2} and Eqn.\,\eqref{eq:pr3} into
Eqn.\,\eqref{eq:pr1}, we have 
\begin{equation}
	\begin{aligned}
&\frac{1}{\Delta t}\|\phi_h^{n+1}-\phi_h^{n}\|^2 + \frac{\xi^2}{2} \|\psi_h^{n+1}\|^2
 -\frac{\xi^2}{2}\|\psi_h^{n}\|^2 {}\\
&+ \frac{\xi^2}{2}\|\psi_h^{n+1}-\psi_h^{n}\|^2 
+(s_h^{n+1})^2-(s_h^{n})^2+(s_h^{n+1}-s_h^{n})^2 =0.
	\end{aligned}
\label{eq:pr4}
\end{equation}
Then the discretized energy dissipative property is satisfied, i.e.
\begin{equation}\label{eq:pr5}
	\begin{aligned}
&\tilde{\mathcal{E}}(\psi_h^{n+1},s_h^{n+1})-\tilde{\mathcal{E}}(\psi_h^{n},s_h^{n}){}\\
&= -\big\{\frac{1}{\Delta t}\|\phi_h^{n+1}-\phi_h^{n}\|^2
 + \frac{\xi^2}{2}\|\psi_h^{n+1}-\psi_h^{n}\|^2 +(s_h^{n+1}-s_h^{n})^2 \big\}{}\\
& \le 0.
	\end{aligned}
\end{equation}

The energy stability of Eqns.\,\eqref{eq:LB-AC-RFOSAV-LFEM} is derived for the modified free energy $\tilde{\mathcal{E}}(\psi,s)$, not for the original one $\mathcal{E}(\phi)$, owing to the introduction of a function $\psi$ and a scaler variable $s$ in the constructing processes of Eqns.\,\eqref{eq:LB-AC-RFOSAV-LFEM}. The modified free energy plays an important role in the error analysis of the discrete scheme.

\subsection{Error estimate}
We first give some lemmas below.

\textbf{Lemma 3.1}\quad  Assume $\phi^0:=\phi(0) \in H^2$. Let $\phi(t)$ be solutions of
Eqns.\,\eqref{eq:LB-AC-SAV}. There exists a constant $C$ depending only on $\Omega$ and $\phi^0$ such that, 
\begin{equation}\label{eq:phit}
\|\phi(t)\|_{H^2} \le C.
\end{equation}

\textbf{Proof:}\quad
It is easy to know that Eqns.\,\eqref{eq:LB-AC-SAV} are unconditionally energy stability with modified energy \eqref{eq:pr5}. There exists a constant $C$ depending only on $\Omega$ and $\phi^0$ such that $\forall\,n$, 
\begin{equation}\label{eq:ng}
\frac{\xi^2}{2}\|\psi(t^{n})\|^2 + |s(t^{n})|^2 \le C.
\end{equation}
We take the $L^2$ inner product of Eqn.\,\eqref{eq:LB-AC-SAV1} with $\phi$, then for all $n$,
\begin{equation}\label{eq:eepr1}
\big(\phi_t(t^{n}), \phi(t^{n})\big)+\xi^2\big((\nabla^2+1) \psi(t^{n}),\phi(t^{n})\big) +\big(u(\phi(t^{n-1})) s(t^{n}),\phi(t^{n})\big) = 0.
\end{equation}
Integrating Eqn.\,\eqref{eq:LB-AC-SAV2} with $\psi(t^{n})$,  we obtain
\begin{equation}\label{eq:eepr2}
\big(\psi(t^{n}), (\nabla^2+1)\phi(t^{n}) \big)=\|\psi(t^{n})\|^2.
\end{equation}
Combining Eqn.\,\eqref{eq:eepr1} and Eqn.\,\eqref{eq:eepr2}, and using H{\"o}lder's inequality and Young's inequality, we have
\begin{equation}\label{eq:eepr3}
\begin{aligned}
\frac{1}{2}\frac{d\|\phi(t^{n})\|^2}{dt}+\|\psi(t^{n})\|^2&= -\big(u(\phi(t^{n-1})) s(t^{n}),\phi(t^{n})\big) \\
&\le |s(t^{n})|\|u(\phi(t^{n-1}))\|\|\phi(t^{n})\|.
\end{aligned}
\end{equation}
Thanks to Eqn.\,\eqref{eq:ng}, then
\begin{equation}\label{eq:eepr31}
\begin{aligned}
\frac{1}{2}\frac{d\|\phi(t^{n})\|}{dt}&\le C\|u(\phi(t^{n-1}))\|.
\end{aligned}
\end{equation}
By Minkowski's inequality and Sobolve embedding, and choosing the appropriate $D_0$ such that $\mathcal{E}_1(\phi(t^{n-1}))+D_0\ge 1$ we have
\begin{equation}\label{eq:eepr4}
\begin{aligned}
\|\phi(t^{n})\| &\le C\Delta t\|u(\phi(t^{n-1}))\|+\|\phi(t^{n-1})\|\\
&= C\Delta t\|\frac{\mathcal{N}'(\phi(t^{n-1}))}{\sqrt{\mathcal{E}_1(\phi(t^{n-1}))+D_0}}\| +\|\phi(t^{n-1})\|\\
&\le C\Delta t\|\mathcal{N}'(\phi(t^{n-1}))\|+\|\phi(t^{n-1})\|\\
&= C\Delta t\|\alpha\phi(t^{n-1})-\frac{\gamma}{2}(\phi(t^{n-1}))^2+\frac{1}{6}(\phi(t^{n-1}))^3\|+\|\phi(t^{n-1})\|\\
&\le C\Delta t(\alpha\|\phi(t^{n-1})\|+\frac{\gamma}{2}\|(\phi(t^{n-1}))^2\|+\frac{1}{6}\|(\phi(t^{n-1}))^3\|)+\|\phi(t^{n-1})\|\\
&\le C\Delta t(\alpha\|\phi(t^{n-1})\|+c_1\frac{\gamma}{2}\|\phi(t^{n-1})\|_{L^6}^2+\frac{1}{6}\|\phi(t^{n-1})\|_{L^6}^3)+ \|\phi(t^{n-1})\|\\
&\le C\Delta t(\alpha\|\phi(t^{n-1})\|+c_1\frac{\gamma}{2}\|\phi(t^{n-1})\|_{H^1}^2+\frac{1}{6}\|\phi(t^{n-1})\|_{H^1}^3)+\|\phi(t^{n-1})\|.
\end{aligned}
\end{equation}
Next we shall prove $\|\phi(t^{n})\| \le C$ using mathematical induction.\\
When $n=0$,  $\|\phi(t^{0})\|=\|\phi^0\| \le C$.\\
If $n=k$, $\|\phi(t^{k})\| \le C$, by Minkowski's inequality and H{\"o}lder's inequality, we have the results as follow:
\begin{equation}\label{eq:eepr5}
\begin{aligned}
&\|\nabla^2 \phi(t^{k})\| = \|(\nabla^2+1) \phi(t^{k})- \phi(t^{k})\| \le \|(\nabla^2+1) \phi(t^{k})\|+ \|\phi(t^{k})\| \le C,\\
&\|\nabla \phi(t^{k})\|^2 = -\big(\nabla^2 \phi(t^{k}), \phi(t^{k}) \big) \le \|\nabla^2 \phi(t^{k})\| \|\phi(t^{k})\| \le C.
\end{aligned}
\end{equation}
Note that $\|\phi(t^{k})\| \le C$ and $\|\phi(t^{k})\|_{H^1} \le C$, using Eqns.\,\eqref{eq:eepr4}, we deduce that $\|\phi(t^{k+1})\|^2 \le C$.\\
Thus, $\|\phi^n\| \le C$. In addition to Eqns.\,\eqref{eq:eepr5}, we can deduce Eqn.\,\eqref{eq:phin}.
\\

For Eqns.\,\eqref{eq:LB-AC-RFOSAV-LFEM}, we can obtain Lemma 3.2,
as its proofs essentially the same as for Lemma 3.1.

\textbf{Lemma 3.2}\quad Assume $\phi^0 \in H^2$. Let $\phi_h^n$ be solutions of Eqns.\,\eqref{eq:LB-AC-RFOSAV-LFEM}. Then for all $h$ and $n < T/{\Delta t}$, we have 
\begin{equation}\label{eq:phin}
\|\phi_h^n\|_{H^2} \le C.
\end{equation}
where $\|\phi_h^n\|_{H^2}^2  = \|\phi_h^n\|^2 + \|\nabla
\phi_h^n\|^2 +\|\Delta_h \phi_h^n\|^2$ and $\Delta_h: H^1 \to
V_h$ is the discrete Laplace operator.

Here, we shall derive error estimates of Eqns.\,\eqref{eq:LB-AC-RFOSAV-LFEM}. 
Denote $e_{\phi,h}^{n}=\phi_h^{n}-\phi(t^n)$, $e_{\psi,h}^{n}=\psi_h^{n}-\psi(t^n)$, $e_{s,h}^{n}=s_h^{n}-s(t^n)$, where $s(t^n)=\sqrt{\mathcal{E}_1(\phi(t^n))+D_0}$.

Let $R_h: H^1(\Omega) \to V_h$ be the standard elliptic (Ritz) projection operator, satisfying 
\begin{equation*}
(\nabla(R_h w(t) - w(t)),\nabla v_h) = 0,\quad \forall v_h \in V_h,
\end{equation*}
and $\theta_{w,h}^n:=w_h^{n}-R_h w(t^n),\quad \rho_{v,h}^n:=R_h w(t^n)-w(t^n),\quad D_t w_h^{n+1} := \frac{w_h^{n+1}-w_h^{n}}{\Delta t}$.

\textbf{Lemma 3.3}\cite{thomee1984galerkin,elliott1989second}\quad If $w$
are sufficiently smooth,  there exist a positive constant $C$ 
independent of $t \in [0,T]$, such that 
\begin{align}
&\|\rho(t)\|+h\|\nabla \rho(t)\| \le Ch^2\|w\|_{H^2},\\
&\|\rho_t(t)\|+h\|\nabla \rho_t(t)\| \le Ch^2\|w_t\|_{H^2}.
\end{align}

\textbf{Theorem 3.2}\quad
Let $\phi$ and $\phi_h^n$ be solutions of Eqns.\,\eqref{eq:LB-AC-SAV} and Eqns.\,\eqref{eq:LB-AC-RFOSAV-LFEM}, respectively. Assume $\phi^0 \in H^2$. In addition, we assume that
$$\phi \in  C(0,T;H^2),\quad \phi_t \in  L^{\infty}(0,T;L^2) \cap L^2(0,T;L^2) \cap C(0,T;H^2),\quad \phi_{tt} \in  L^2(0,T;L^2),$$
and 
$$\psi \in  C(0,T;H^2),\quad  \psi_t \in C(0,T;H^2)$$
For $0 < n < T/\Delta t$, such that
\begin{equation}\label{eq:eest}
\|e_{\phi,h}^n\|^2+\frac{\xi^2}{2} \| e_{\psi,h}^{n}\|^2 +(e_{s,h}^n)^2 \le C(K_1 \Delta t^2+ K_2 h^4).
\end{equation} 
where 
\begin{align}
K_1 &= \|\phi_t\|_{L^2(0,T;L^2)}^4+\|\phi_{tt}\|_{L^2(0,T;L^2)}^2,\nonumber\\
K_2 &= \|\phi\|_{C(0,T;H^2)}^2+\|\phi_t\|_{C(0,T;H^2)}^2+\|\psi\|_{C(0,T;H^2)}^2+\|\psi_t\|_{C(0,T;H^2)}^2.\nonumber
\end{align}

\textbf{Proof:}\quad
By Lemma 3.1 and  Lemma 3.2, we know that 
$$\|\phi(t)\|_{H^2},\|\phi_h^n\|_{H^2}  \le C.$$
Note that $H^2 \subseteq L^{\infty}$. Therefore, we can find a constant $C$ such that
\begin{align}\label{eq:up}
|\mathcal{N}(\phi(t))|,|\mathcal{N}'(\phi(t))|,|\mathcal{N}''(\phi(t))| \le C,\nonumber\\
|\mathcal{N}(\phi_h^n)|,|\mathcal{N}'(\phi_h^n)|,|\mathcal{N}''(\phi_h^n)| \le C.
\end{align}
Subtracting the variational formulation of Eqn.\,\eqref{eq:LB-AC-SAV1} from Eqn.\,\eqref{eq:LB-AC-RFOSAV-LFEM1}, we have
\begin{align}\label{eq:err1a}
&\big(D_t e_{\phi,h}^{n+1}+\xi^2 e_{\psi,h}^{n+1},v_h\big)
-\xi^2\big(\nabla e_{\psi,h}^{n+1},\nabla v_h\big) \nonumber\\
&= - \big(u(\phi_h^{n}) s_h^{n+1}-u(\phi(t^{n})) s(t^{n+1}),v_h\big)  - \big(D_t \phi(t^{n+1})-\phi_t(t^{n+1}) ,v_h\big).
\end{align}
Then 
\begin{align}\label{eq:err1}
&\big(D_t \theta_{\phi,h}^{n+1}+\xi^2 \theta_{\psi,h}^{n+1},v_h\big)
-\xi^2\big(\nabla \theta_{\psi,h}^{n+1},\nabla v_h\big) \nonumber\\
&= - \big(u(\phi_h^{n}) s_h^{n+1}-u(\phi(t^{n})) s(t^{n+1}),v_h\big)  - \big(D_t \phi(t^{n+1})-\phi_t(t^{n+1}),v_h\big)\nonumber\\
&\quad - \big(D_t \rho_{\phi,h}^{n+1}+\xi^2 \rho_{\psi,h}^{n+1},v_h\big)
-\xi^2\big(\nabla \rho_{\psi,h}^{n+1},\nabla v_h\big)\nonumber\\
&= - \big(u(\phi_h^{n}) s_h^{n+1}-u(\phi(t^{n})) s(t^{n+1}),v_h\big)  - \big(D_t \phi(t^{n+1})-\phi_t(t^{n+1}) ,v_h\big)\nonumber\\
&\quad - \big(D_t \rho_{\phi,h}^{n+1}+\xi^2\rho_{\psi,h}^{n+1},v_h\big).
\end{align}
Subtracting the variational formulation of Eqn.\,\eqref{eq:LB-AC-SAV2} from Eqn.\,\eqref{eq:LB-AC-RFOSAV-LFEM2}, we get
\begin{align}\label{eq:err2a}
\big(e_{\psi,h}^{n+1},v_h\big)=\big(e_{\phi,h}^{n+1},v_h\big)-\big(\nabla e_{\phi,h}^{n+1},\nabla v_h\big).
\end{align}
Then
\begin{align}\label{eq:err2}
&\big(\theta_{\phi,h}^{n+1},v_h\big)-\big(\nabla \theta_{\phi,h}^{n+1},\nabla v_h\big)\nonumber\\
&=\big(\theta_{\psi,h}^{n+1},v_h\big)+\big(\rho_{\psi,h}^{n+1},v_h\big)-\big(\rho_{\phi,h}^{n+1},v_h\big)+\big(\nabla\rho_{\phi,h}^{n+1},\nabla v_h\big)\nonumber\\
&=\big(\theta_{\psi,h}^{n+1},v_h\big)+\big(\rho_{\psi,h}^{n+1}-\rho_{\phi,h}^{n+1},v_h\big).
\end{align}
Taking $v_h = \theta_{\phi,h}^{n+1}$ in Eqn.\,\eqref{eq:err1} and $v_h=\theta_{\psi,h}^{n+1}$ in 
Eqn.\,\eqref{eq:err2}, we obtain
\begin{align}\label{eq:2err}
&\frac{1}{2}D_t \|\theta_{\phi,h}^{n+1}\|^2+\frac{\|\theta_{\phi,h}^{n+1}-\theta_{\phi,h}^{n}\|^2}{2\Delta t}+\xi^2 \|\theta_{\psi,h}^{n+1}\|^2 \nonumber\\
&= - \big(u(\phi_h^{n}) s_h^{n+1}-u(\phi(t^{n})) s(t^{n+1}),\theta_{\phi,h}^{n+1}\big) \nonumber\\
&\quad - \big(D_t \phi(t^{n+1})-\phi_t(t^{n+1}) ,\theta_{\phi,h}^{n+1}\big)\nonumber\\
&\quad - \big(D_t \rho_{\phi,h}^{n+1}+\xi^2\rho_{\psi,h}^{n+1},\theta_{\phi,h}^{n+1}\big)\nonumber\\
&\quad - \xi^2 \big(\rho_{\psi,h}^{n+1}-\rho_{\phi,h}^{n+1},\theta_{\psi,h}^{n+1}\big)\nonumber\\
&= \big(u(\phi_h^{n})(s_h^{n+1}-s(t^{n+1})),\theta_{\phi,h}^{n+1}\big)\nonumber\\
&\quad+s(t^{n+1})\big(u(\phi_h^{n}) -u(\phi(t^{n})), \theta_{\phi,h}^{n+1}\big)\nonumber\\
&\quad - \big(D_t \phi(t^{n+1})-\phi_t(t^{n+1}) ,\theta_{\phi,h}^{n+1}\big)\nonumber\\
&\quad - \big(D_t \rho_{\phi,h}^{n+1}+\xi^2\rho_{\psi,h}^{n+1},\theta_{\phi,h}^{n+1}\big)\nonumber\\
&\quad - \xi^2 \big(\rho_{\psi,h}^{n+1}-\rho_{\phi,h}^{n+1},\theta_{\psi,h}^{n+1}\big)\nonumber\\
&:=I_1+I_2+\cdots+I_5.
\end{align}
Using H{\"o}lder's inequality, we get
\begin{align}\label{eq:I1}
I_1&=\big(u(\phi_h^{n})e_{s,h}^{n+1},\theta_{\phi,h}^{n+1}\big)\le C(e_{s,h}^{n+1})^2+\|\theta_{\phi,h}^{n+1}\|^2.
\end{align}
Due to
\begin{align*}
&u(\phi_h^n)-u(\phi(t^{n}))\nonumber\\
&=\frac{\mathcal{N}'(\phi_h^n)}{\sqrt{\mathcal{E}_1(\phi_h^n)+D_0}}-\frac{\mathcal{N}'(\phi(t^n))}{\sqrt{\mathcal{E}_1(\phi(t^n))+D_0}}\nonumber\\
&=\frac{\mathcal{N}'(\phi_h^n)}{\sqrt{\mathcal{E}_1(\phi_h^n)+D_0}}-\frac{\mathcal{N}'(\phi_h^n)}{\sqrt{\mathcal{E}_1(\phi(t^n))+D_0}}+\frac{\mathcal{N}'(\phi_h^n)}{\sqrt{\mathcal{E}_1(\phi(t^n))+D_0}}-\frac{\mathcal{N}'(\phi(t^n))}{\sqrt{\mathcal{E}_1(\phi(t^n))+D_0}}\nonumber\\
&=\frac{\mathcal{N}'(\phi_h^n)(\mathcal{E}_1(\phi(t^n))-\mathcal{E}_1(\phi_h^{n}))}{\sqrt{\mathcal{E}_1(\phi_h^n)+D_0}\sqrt{\mathcal{E}_1(\phi(t^n))+D_0}(\sqrt{\mathcal{E}_1(\phi_h^n)+D_0}+\sqrt{\mathcal{E}_1(\phi(t^n))+D_0})}\nonumber\\
&\quad + \frac{\mathcal{N}'(\phi_h^n)-\mathcal{N}'(\phi(t^n))}{\sqrt{\mathcal{E}_1(\phi(t^n))+D_0}}.
\end{align*}
Using inequalities \eqref{eq:up}, we derive
\begin{align*}
\|u(\phi_h^n)-u(\phi(t^{n}))\| 
&\le C(\|\mathcal{E}_1(\phi(t^n))-\mathcal{E}_1(\phi_h^{n}))\|+\|\mathcal{N}'(\phi_h^n)-\mathcal{N}'(\phi(t^n))\|)
\nonumber\\
&\le C (\|\mathcal{N}(\eta_1)(\phi_h^{n}-\phi(t^{n}))\|+\|\mathcal{N}''(\eta_2)(\phi_h^{n}-\phi(t^{n}))\|)\nonumber\\
&\le C\|e_{\phi,h}^{n}\|\nonumber\\
& \le C(\|\rho_{\phi,h}^{n}\|+\|\theta_{\phi,h}^{n}\|)\nonumber\\
&\le Ch^2\|\phi\|_{H^2}+C\|\theta_{\phi,h}^{n}\|,
\end{align*}
where $\eta_1,\eta_2$ lies between $\phi_h^{n}$ and $\phi(t^{n})$. 
Note that $s(t^{n+1})<C$, we find 
\begin{align}\label{eq:I2}
I_2 &\le C \|u(\phi_h^n)-u(\phi(t^{n}))\|\|\theta_{\phi,h}^{n+1}\|\nonumber\\
&\le C \|u(\phi_h^n)-u(\phi(t^{n}))\|^2+\|\theta_{\phi,h}^{n+1}\|^2\nonumber\\
&\le Ch^4\|\phi\|_{H^2}^2+C\|\theta_{\phi,h}^{n}\|^2+\|\theta_{\phi,h}^{n+1}\|^2.
\end{align}
For $I_3$, we have
\begin{align}
I_3&\le \|D_t \phi(t^{n+1})-\phi_t(t^{n+1})\|\|\theta_{\phi,h}^{n+1}\|\nonumber\\
&\le \frac{1}{4}\|D_t \phi(t^{n+1})-\phi_t(t^{n+1})\|^2+\|\theta_{\phi,h}^{n+1}\|^2\nonumber\\
&\le \frac{\Delta t}{4}\int_{t^n}^{t^{n+1}}\|\phi_{tt}(r)\|^2\mathrm d r+\|\theta_{\phi,h}^{n+1}\|^2.\label{eq:I3}
\end{align}
$I_4$ and $I_5$ can be estimated as follow:
\begin{align}
I_4&\le \|D_t \rho_{\phi,h}^{n+1}+\xi^2\rho_{\psi,h}^{n+1}\|\|\theta_{\phi,h}^{n+1}\| \nonumber\\
&\le  \frac{1}{4}\|D_t \rho_{\phi,h}^{n+1}+\xi^2\rho_{\psi,h}^{n+1}\|^2+\|\theta_{\phi,h}^{n+1}\|^2 \nonumber\\
&\le \frac{1}{4}\|D_t \rho_{\phi,h}^{n+1}\|^2+ \frac{\xi^2}{4}\|\rho_{\psi,h}^{n+1}\|^2+\|\theta_{\phi,h}^{n+1}\|^2 \nonumber\\
&\le Ch^4(\|\phi_t\|_{H^2}^2+\|\psi\|_{H^2}^2)+\|\theta_{\phi,h}^{n+1}\|^2,\label{eq:I4}\\
I_5&\le \xi^2\|\rho_{\psi,h}^{n+1}-\rho_{\phi,h}^{n+1}\|\|\theta_{\psi,h}^{n+1}\| \nonumber\\
&\le \frac{\xi^2}{2}\|\rho_{\psi,h}^{n+1}-\rho_{\phi,h}^{n+1}\|^2+\xi^2\|\theta_{\psi,h}^{n+1}\|^2 \nonumber\\
&\le \frac{\xi^2}{2}(\|\rho_{\psi,h}^{n+1}\|^2+\|\rho_{\phi,h}^{n+1}\|^2)+\frac{\xi^2}{2}\|\theta_{\psi,h}^{n+1}\|^2\nonumber\\
&\le Ch^4(\|\psi\|_{H^2}^2+\|\phi\|_{H^2}^2)+\frac{\xi^2}{2}\|\theta_{\psi,h}^{n+1}\|^2.\label{eq:I5}
\end{align}
Using inequalities \eqref{eq:I1}-\eqref{eq:I5}, Eqn.\,\eqref{eq:2err} can be estimated as
\begin{align}\label{eq:re1}
&\|\theta_{\phi,h}^{n+1}\|^2-\|\theta_{\phi,h}^{n}\|^2+\|\theta_{\phi,h}^{n+1}-\theta_{\phi,h}^{n}\|^2+\xi^2 \Delta t \|\theta_{\psi,h}^{n+1}\|^2
\le(C\Delta t^2\int_{t^n}^{t^{n+1}}\|\phi_{tt}(r)\|^2\mathrm d r\nonumber\\
&\quad+C\Delta th^4(\|\phi\|_{H^2}^2+\|\phi_t\|_{H^2}^2+\|\psi\|_{H^2}^2))
+ C\Delta t(\|\theta_{\phi,h}^{n}\|^2+\|\theta_{\phi,h}^{n+1}\|^2+(e_{s,h}^{n+1})^2).
\end{align}
Taking $v_h = \theta_{\phi,h}^{n+1}-\theta_{\phi,h}^{n}$ in Eqn.\,\eqref{eq:err1}, we obtain
\begin{align}\label{eq:err2a}
&\frac{\|\theta_{\phi,h}^{n+1}-\theta_{\phi,h}^{n}\|^2}{\Delta t}+\xi^2 \big(\theta_{\psi,h}^{n+1},\theta_{\phi,h}^{n+1}-\theta_{\phi,h}^{n}\big)-\xi^2 \big(\nabla \theta_{\psi,h}^{n+1},\nabla (\theta_{\phi,h}^{n+1}-\theta_{\phi,h}^{n})\big)  \nonumber\\
&= - \big(u(\phi_h^{n}) s_h^{n+1}-u(\phi(t^{n})) s(t^{n+1}),\theta_{\phi,h}^{n+1}-\theta_{\phi,h}^{n}\big) \nonumber\\
&\quad - \big(D_t \phi(t^{n+1})-\phi_t(t^{n+1}) ,\theta_{\phi,h}^{n+1}-\theta_{\phi,h}^{n}\big)\nonumber\\
&\quad - \big(D_t \rho_{\phi,h}^{n+1}+\xi^2\rho_{\psi,h}^{n+1},\theta_{\phi,h}^{n+1}-\theta_{\phi,h}^{n}\big).
\end{align}
Taking $v_h=\theta_{\psi,h}^{n+1}$ to Eq.\,\eqref{eq:err2}, we know
\begin{align}\label{eq:err2b}
&\big(\theta_{\psi,h}^{n+1},\theta_{\phi,h}^{n+1}-\theta_{\phi,h}^{n}\big)- \big(\nabla \theta_{\psi,h}^{n+1},\nabla (\theta_{\phi,h}^{n+1}-\theta_{\phi,h}^{n})\big)\nonumber\\
&=\big(\theta_{\psi,h}^{n+1},\theta_{\phi,h}^{n+1}\big)- \big(\nabla \theta_{\psi,h}^{n+1},\nabla \theta_{\phi,h}^{n+1}\big)-[\big(\theta_{\psi,h}^{n+1},\theta_{\phi,h}^{n}\big)- \big(\nabla \theta_{\psi,h}^{n+1},\nabla \theta_{\phi,h}^{n}\big)]\nonumber\\
&= \big(\theta_{\psi,h}^{n+1}, \theta_{\psi,h}^{n+1}\big)+\big(\rho_{\psi,h}^{n+1}-\rho_{\phi,h}^{n+1},\theta_{\psi,h}^{n+1}\big)
- \big(\theta_{\psi,h}^{n},\theta_{\psi,h}^{n+1}\big)-\big(\rho_{\psi,h}^{n}-\rho_{\phi,h}^{n},\theta_{\psi,h}^{n+1}\big)\nonumber\\
&=\big( \theta_{\psi,h}^{n+1},\theta_{\psi,h}^{n+1}-\theta_{\psi,h}^{n}\big)+\big( \theta_{\psi,h}^{n+1},\rho_{\psi,h}^{n+1}-\rho_{\psi,h}^{n}\big)-\big( \theta_{\psi,h}^{n+1},\rho_{\phi,h}^{n+1}-\rho_{\phi,h}^{n}\big).
\end{align}
Substituting Eqn.\,\eqref{eq:err2b} into Eqn.\,\eqref{eq:err2a}, we have
\begin{align}\label{eq:2err2}
&\frac{\|\theta_{\phi,h}^{n+1}-\theta_{\phi,h}^{n}\|^2}{\Delta t}+\frac{\xi^2}{2} (\|\theta_{\psi,h}^{n+1}\|^2-\|\theta_{\psi,h}^{n}\|^2)+\frac{\xi^2}{2}\|\theta_{\psi,h}^{n+1}-\theta_{\psi,h}^{n}\|^2 \nonumber\\
&= - \big(u(\phi_h^{n})e_{s,h}^{n+1},\theta_{\phi,h}^{n+1}-\theta_{\phi,h}^{n}\big) \nonumber\\
&\quad- \big((u(\phi_h^{n})-u(\phi(t^{n}))) s(t^{n+1}),\theta_{\phi,h}^{n+1}-\theta_{\phi,h}^{n}\big) \nonumber\\
&\quad - \big(D_t \phi(t^{n+1})-\phi_t(t^{n+1}) ,\theta_{\phi,h}^{n+1}-\theta_{\phi,h}^{n}\big)\nonumber\\
&\quad - \big(D_t \rho_{\phi,h}^{n+1}+\xi^2\rho_{\psi,h}^{n+1},\theta_{\phi,h}^{n+1}-\theta_{\phi,h}^{n}\big)\nonumber\\
&\quad -\big( \theta_{\psi,h}^{n+1},\rho_{\psi,h}^{n+1}-\rho_{\psi,h}^{n}\big) +\big( \theta_{\psi,h}^{n+1},\rho_{\phi,h}^{n+1}-\rho_{\phi,h}^{n}\big).
\end{align}
Subtracting Eqn.\,\eqref{eq:LB-AC-SAV3} from Eqn.\,\eqref{eq:LB-AC-RFOSAV-LFEM3}, we have
\begin{align}\label{eq:s1}
 &e_{s,h}^{n+1}-e_{s,h}^{n}\nonumber\\
 &= \frac{1}{2}[\big(u(\phi_h^{n}), \phi_h^{n+1}-\phi_h^{n}\big)-\big(u(\phi(t^{n})), \Delta t \phi_t(t^{n+1}) \big)]\nonumber\\
&\quad -(s(t^{n+1})-s(t^{n})-\Delta t s_t(t^{n+1}))\nonumber\\
&= \frac{1}{2}\big(u(\phi_h^{n}),\phi_h^{n+1}-\phi_h^{n}-\Delta t\phi_t(t^{n+1})\big)+\frac{1}{2}\big(u(\phi_h^{n})-u(\phi(t^{n})), \Delta t \phi_t(t^{n+1}) \big)\nonumber\\
&\quad -(s(t^{n+1})-s(t^{n})-\Delta t s_t(t^{n+1}))\nonumber\\
&= \frac{1}{2}\big(u(\phi_h^{n}),e_{\phi,h}^{n+1}-e_{\phi,h}^n\big)+\frac{1}{2}\big(u(\phi_h^{n}),\phi(t^{n+1})-\phi(t^{n})-\Delta t\phi_t(t^{n+1})\big)\nonumber\\
&\quad+\frac{1}{2}\big(u(\phi_h^{n})-u(\phi(t^{n})), \Delta t\phi_t(t^{n+1}) \big)-(s(t^{n+1})-s(t^{n})-\Delta t s_t(t^{n+1}))\nonumber\\
&= \frac{1}{2}\big(u(\phi_h^{n}),\theta_{\phi,h}^{n+1}-\theta_{\phi,h}^n\big)+\frac{1}{2}\big(u(\phi_h^{n}),\rho_{\phi,h}^{n+1}-\rho_{\phi,h}^n\big)\nonumber\\
&\quad+\frac{1}{2}\big(u(\phi_h^{n}),\phi(t^{n+1})-\phi(t^{n})-\Delta t\phi_t(t^{n+1})\big)\nonumber\\
&\quad+\frac{1}{2}\big(u(\phi_h^{n})-u(\phi(t^{n})), \Delta t\phi_t(t^{n+1}) \big)-(s(t^{n+1})-s(t^{n})-\Delta t s_t(t^{n+1})).
\end{align}
Multiplying Eqn.\,\eqref{eq:s1} by $2e_{s,h}^{n+1}$, then
\begin{align}\label{eq:s2}
&(e_{s,h}^{n+1})^2-(e_{s,h}^{n})^2+(e_{s,h}^{n+1}-e_{s,h}^{n})^2 \nonumber\\
&= e_{s,h}^{n+1}\big(u(\phi_h^{n}),\theta_{\phi,h}^{n+1}-\theta_{\phi,h}^n\big)+e_{s,h}^{n+1}\big(u(\phi_h^{n}),\rho_{\phi,h}^{n+1}-\rho_{\phi,h}^n\big)\nonumber\\
&\quad+e_{s,h}^{n+1}\big(u(\phi_h^{n}),\phi(t^{n+1})-\phi(t^{n})-\Delta t\phi_t(t^{n+1})\big)\nonumber\\
&\quad+e_{s,h}^{n+1}\big(u(\phi_h^{n})-u(\phi(t^{n})), \Delta t\phi_t(t^{n+1}) \big)\nonumber\\
&\quad-2e_{s,h}^{n+1}(s(t^{n+1})-s(t^{n})-\Delta t s_t(t^{n+1})).
\end{align}
Combining Eqn.\,\eqref{eq:s2} and Eqn.\,\eqref{eq:2err2}, we obtain
\begin{align}\label{eq:2erra}
&\frac{\|\theta_{\phi,h}^{n+1}-\theta_{\phi,h}^{n}\|^2}{\Delta t}+\frac{\xi^2}{2} (\|\theta_{\psi,h}^{n+1}\|^2-\|\theta_{\psi,h}^{n}\|^2)+\frac{\xi^2}{2}\|\theta_{\psi,h}^{n+1}-\theta_{\psi,h}^{n}\|^2\nonumber\\
&\quad + (e_{s,h}^{n+1})^2-(e_{s,h}^{n})^2+(e_{s,h}^{n+1}-e_{s,h}^{n})^2 \nonumber\\
&= - \big((u(\phi_h^{n})-u(\phi(t^{n}))) s(t^{n+1}),\theta_{\phi,h}^{n+1}-\theta_{\phi,h}^{n}\big) \nonumber\\
&\quad - \big(D_t \phi(t^{n+1})-\phi_t(t^{n+1}) ,\theta_{\phi,h}^{n+1}-\theta_{\phi,h}^{n}\big)\nonumber\\
&\quad - \big(D_t \rho_{\phi,h}^{n+1}+\xi^2\rho_{\psi,h}^{n+1},\theta_{\phi,h}^{n+1}-\theta_{\phi,h}^{n}\big)\nonumber\\
&\quad -\big( \theta_{\psi,h}^{n+1},\rho_{\psi,h}^{n+1}-\rho_{\psi,h}^{n}\big) \nonumber\\
&\quad+\big( \theta_{\psi,h}^{n+1},\rho_{\phi,h}^{n+1}-\rho_{\phi,h}^{n}\big)\nonumber\\
&\quad+e_{s,h}^{n+1}\big(u(\phi_h^{n}),\rho_{\phi,h}^{n+1}-\rho_{\phi,h}^n\big)\nonumber\\
&\quad+e_{s,h}^{n+1}\big(u(\phi_h^{n}),\phi(t^{n+1})-\phi(t^{n})-\Delta t\phi_t(t^{n+1})\big)\nonumber\\
&\quad+e_{s,h}^{n+1}\big(u(\phi_h^{n})-u(\phi(t^{n})), \Delta t\phi_t(t^{n+1}) \big)\nonumber\\
&\quad-2e_{s,h}^{n+1}(s(t^{n+1})-s(t^{n})-\Delta t s_t(t^{n+1}))\nonumber\\
&:=J_1+J_2+\cdots+J_9.
\end{align}
For $J_1$ and $J_2$, we have
\begin{align}\label{eq:J12}
J_1&\le |s(t^{n+1})|\|u(\phi_h^{n})-u(\phi(t^{n}))\|\|\theta_{\phi,h}^{n+1}-\theta_{\phi,h}^{n}\|\nonumber\\
&\le C\Delta t \|u(\phi_h^n)-u(\phi(t^{n}))\|^2+\frac{1}{4\Delta t}\|\theta_{\phi,h}^{n+1}-\theta_{\phi,h}^{n}\|^2\nonumber\\
&\le C\Delta t (\|\theta_{\phi}^{n}\|^2+ h^4\|\phi\|_{H^2}^2)+\frac{1}{4\Delta t}\|\theta_{\phi,h}^{n+1}-\theta_{\phi,h}^{n}\|^2,\\
J_2&\le \|D_t \phi(t^{n+1})-\phi_t(t^{n+1})\|\|\theta_{\phi,h}^{n+1}-\theta_{\phi,h}^{n}\|\nonumber\\&\le \Delta t\|D_t \phi(t^{n+1})-\phi_t(t^{n+1})\|^2+ \frac{1}{4\Delta t}\|\theta_{\phi,h}^{n+1}-\theta_{\phi,h}^{n}\|^2\nonumber\\
&\le  \Delta t^2 \int_{t^n}^{t^{n+1}}\|\phi_{tt}(r)\|^2\mathrm d r+\frac{1}{4\Delta t}\|\theta_{\phi,h}^{n+1}-\theta_{\phi,h}^{n}\|^2.
\end{align}
For $J_3$ and $J_4$, we have
\begin{align}\label{eq:J34}
J_3 &\le \|D_t \rho_{\phi,h}^{n+1}+\xi^2 \rho_{\psi,h}^{n+1}\|\|\theta_{\phi,h}^{n+1}-\theta_{\phi,h}^{n}\|\nonumber\\
&\le \Delta t(\|D_t \rho_{\phi,h}^{n+1}\|^2+\xi^2 \|\rho_{\psi,h}^{n+1}\|^2)+\frac{1}{4\Delta t}\|\theta_{\phi,h}^{n+1}-\theta_{\phi,h}^{n}\|^2\nonumber\\
&\le C\Delta t h^4(\|\phi_t\|_{H^2}^2+\|\psi\|_{H^2}^2)+\frac{1}{4\Delta t}\|\theta_{\phi,h}^{n+1}-\theta_{\phi,h}^{n}\|^2,\\
J_4&\le\|\theta_{\psi,h}^{n+1}\|\|\rho_{\psi,h}^{n+1}-\rho_{\psi,h}^{n}\|\nonumber\\
&\le \Delta t\|\theta_{\psi,h}^{n+1}\|^2+\frac{1}{4\Delta t}\|\Delta t D_t \rho_{\psi,h}^{n+1} \|^2\nonumber\\
&\le \Delta t\|\theta_{\psi,h}^{n+1}\|^2+C\Delta t h^4\|\psi_t\|_{H^2}^2.
\end{align}
For $J_5$ and $J_6$, we get
\begin{align}\label{eq:J56}
J_5&\le\|\theta_{\psi,h}^{n+1}\|\|\rho_{\phi,h}^{n+1}-\rho_{\phi,h}^{n}\|\nonumber\\
&\le \Delta t\|\theta_{\psi,h}^{n+1}\|^2+\frac{1}{4\Delta t}\|\Delta t D_t \rho_{\phi,h}^{n+1} \|^2\nonumber\\
&\le \Delta t\|\theta_{\psi,h}^{n+1}\|^2+C\Delta t h^4\|\phi_t\|_{H^2}^2,\\
J_6 &\le \|u(\phi_h^{n})\||e_{s,h}^{n+1}|\|\rho_{\phi,h}^{n+1}-\rho_{\phi,h}^n\|\nonumber\\
 &\le C\Delta t (e_{s,h}^{n+1})^2 + \frac{1}{\Delta t}\|\rho_{\phi,h}^{n+1}-\rho_{\phi,h}^n\|^2\nonumber\\
 &\le C\Delta t (e_{s,h}^{n+1})^2 + \frac{1}{\Delta t} \|\Delta t D_t \rho_{\phi,h}^{n+1}\|^2\nonumber\\
  &\le C\Delta t (e_{s,h}^{n+1})^2 + C \Delta t h^4\|\phi_t\|_{H^2}.
\end{align}
For $J_7,J_8$ ad $J_9$, we have
\begin{align}
J_7&\le |e_{s,h}^{n+1}|\|u(\phi_h^{n})\|\|\phi(t^{n+1})-\phi(t^{n})-\Delta t\phi_t(t^{n+1})\|\nonumber\\
&\le C\Delta t (e_{s,h}^{n+1})^2 + \frac{1}{\Delta t}\|\phi(t^{n+1})-\phi(t^{n})-\Delta t\phi_t(t^{n+1})\|^2\nonumber\\
&\le C\Delta t (e_{s,h}^{n+1})^2 +\Delta t^2 \int_{t^{n}}^{t^{n+1}}\|\phi_{tt}(r)\| \mathrm d r,\\
J_8&\le\Delta t|e_{s,h}^{n+1}|\|u(\phi_h^{n})-u(\phi(t^{n}))\|\|\phi_t(t^{n+1})\|\nonumber\\
&\le C\Delta t\|\phi_t\|_{L^{\infty}(0,T;L^2)}((e_{s,h}^{n+1})^2 + \|u(\phi_h^{n})-u(\phi(t^{n}))\|^2)\nonumber\\
&\le C\Delta t\|\phi_t\|_{L^{\infty}(0,T;L^2)}((e_{s,h}^{n+1})^2 + C\|\theta_{\phi}^{n}\|^2+Ch^4\|\phi\|_{H^2}^2),\\
J_9&\le 2|e_{s,h}^{n+1}|\|s(t^{n+1})-s(t^{n})-\Delta t s_t(t^{n+1})\|\nonumber\\
&\le C\Delta t(e_{s,h}^{n+1})^2 + \frac{1}{\Delta t}\|s(t^{n+1})-s(t^{n})-\Delta t s_t(t^{n+1})\|^2)\nonumber\\
&\le C\Delta t(e_{s,h}^{n+1})^2 + \Delta t^2\int_{t^n}^{t^{n+1}}|s_{tt}(r)|^2\mathrm dr. \label{eq:J9}
\end{align}
In addition, we know
\begin{align}
s_{tt}&=-\frac{1}{4\sqrt{(\mathcal{E}_1(\phi)+D_0)^{3}}}(\int_{\Omega} \mathcal{N}'(\phi)\phi_t\mathrm dr)^2\nonumber\\
&\quad+\frac{1}{2\sqrt{\mathcal{E}_1(\phi)+D_0}}\int_{\Omega}[\mathcal{N}''(\phi)\phi_t^2+\mathcal{N}'(\phi)\phi_{tt}]\mathrm dr,\nonumber\\
\int_{t^n}^{t^{n+1}}|s_{tt}(r)|^2\mathrm dr&\le C\int_{t^n}^{t^{n+1}}\|\phi_t(\tau)\|^4+\|\phi_{tt}(\tau)\|^2\mathrm d\tau. \label{eq:stt1}
\end{align}
Using inequalities \eqref{eq:2erra}-\eqref{eq:J9} and \eqref{eq:stt1}, we obatin
\begin{align}\label{eq:re2}
&\frac{\|\theta_{\phi,h}^{n+1}-\theta_{\phi,h}^{n}\|^2}{4\Delta t}+\frac{\xi^2}{2} (\|\theta_{\psi,h}^{n+1}\|^2-\|\theta_{\psi,h}^{n}\|^2)+\frac{\xi^2}{2}\|\theta_{\psi,h}^{n+1}-\theta_{\psi,h}^{n}\|^2\nonumber\\
&\quad + (e_{s,h}^{n+1})^2-(e_{s,h}^{n})^2+(e_{s,h}^{n+1}-e_{s,h}^{n})^2 \nonumber\\
&\le C(\Delta t^2\int_{t^n}^{t^{n+1}}\|\phi_t(\tau)\|^4+\|\phi_{tt}(\tau)\|^2\mathrm d\tau\nonumber\\ 
&\quad+ \Delta t h^4 (\|\phi_t\|_{H^2}^2+\|\phi\|_{H^2}^2+\|\psi\|_{H^2}^2+\|\psi_t\|_{H^2}^2))\nonumber\\ 
&\quad+C\Delta t(\|\theta_{\psi,h}^{n+1}\|^2+\|\theta_{\phi,h}^{n+1}\|^2+\|\theta_{\phi,h}^{n}\|^2+(e_{s,h}^{n+1})^2).
\end{align}
Combining Eqn.\,\eqref{eq:re1} and Eqn.\,\eqref{eq:re2}, ignoring some nonnegative terms, then
\begin{align}\label{eq:rr}
&\|\theta_{\phi,h}^{n+1}\|^2-\|\theta_{\phi,h}^{n}\|^2+\frac{\xi^2}{2} (\|\theta_{\psi,h}^{n+1}\|^2-\|\theta_{\psi,h}^{n}\|^2) + (e_{s,h}^{n+1})^2-(e_{s,h}^{n})^2\nonumber\\ 
&\le C(\Delta t^2\int_{t^n}^{t^{n+1}}\|\phi_t(\tau)\|^4+\|\phi_{tt}(\tau)\|^2\mathrm d\tau+ \Delta t h^4 (\|\phi\|_{H^2}^2+\|\phi_t\|_{H^2}^2+\|\psi\|_{H^2}^2+\|\psi_t\|_{H^2}^2))\nonumber\\ 
&\quad+C\Delta t(\|\theta_{\psi,h}^{n+1}\|^2+\|\theta_{\phi,h}^{n+1}\|^2+\|\theta_{\phi,h}^{n}\|^2+(e_{s,h}^{n+1})^2).
\end{align}
Adding some nonnegative terms to the right, we find
\begin{align}\label{eq:rr1}
&\|\theta_{\phi,h}^{n+1}\|^2-\|\theta_{\phi,h}^{n}\|^2+\frac{\xi^2}{2} (\|\theta_{\psi,h}^{n+1}\|^2-\|\theta_{\psi,h}^{n}\|^2) + (e_{s,h}^{n+1})^2-(e_{s,h}^{n})^2\nonumber\\ 
&\le C(\Delta t^2\int_{t^n}^{t^{n+1}}\|\phi_t(\tau)\|^4+\|\phi_{tt}(\tau)\|^2\mathrm d\tau 
+ \Delta t h^4 (\|\phi\|_{H^2}^2+\|\phi_t\|_{H^2}^2+\|\psi\|_{H^2}^2+\|\psi_t\|_{H^2}^2))\nonumber\\ 
&\quad+C\Delta t(\|\theta_{\phi,h}^{n+1}\|^2+\|\theta_{\phi,h}^{n}\|^2+
\frac{\xi^2}{2} \|\theta_{\psi,h}^{n+1}\|^2+\frac{\xi^2}{2}\|\theta_{\psi,h}^{n}\|^2+(e_{s,h}^{n+1})^2+(e_{s,h}^{n})^2).
\end{align}
Applying the discrete Gronwall's inequality\cite{shen2018convergence}, we can obtain
\begin{align}\label{eq:rr2}
&\|\theta_{\phi,h}^{n}\|^2+\frac{\xi^2}{2} \|\theta_{\psi,h}^{n}\|^2 + (e_{s,h}^{n})^2\nonumber\\ 
&\le Cexp((1-C\Delta t)^{-1}t^n)(\Delta t^2\int_{0}^{t^n}\|\phi_t(\tau)\|^4+\|\phi_{tt}(\tau)\|^2\mathrm d\tau\nonumber\\ 
&\quad+t^n h^4 (\|\phi\|_{H^2}^2+\|\phi_t\|_{H^2}^2+\|\psi\|_{H^2}^2+\|\psi_t\|_{H^2}^2)).
\end{align}
Combine with Lemma 3.3, thus it is easy to get Eqn.\,\eqref{eq:eest}.

For $H^{-1}$ type dynamic flow, 
\begin{equation}
\label{eq:LB-CH}
\begin{aligned}
\phi_t &= \nabla^2 \varphi,
\\
\varphi &= \xi^2(\nabla^2+1)\psi + u(\phi)s,
\\
\psi &= (\nabla^2+1)\phi,
\\
s_t &= \frac{1}{2}\int_{\Omega} u(\phi)\phi_t \, d\mathbf{r},
\end{aligned}
\end{equation}
its first-order fully discretized scheme is 
\begin{equation}\label{eq:LB-CH-discrete}
\begin{aligned}
\left(\frac{\phi_h^{n+1}-\phi_h^{n}}{\Delta t},v_h\right) &= -\big(\nabla \varphi_h^{n+1},\nabla v_h\big) ,\\
\left( \varphi_h^{n+1},v_h\right) &= \xi^2\big(\psi_h^{n+1},v_h\big)-\xi^2\big(\nabla \psi_h^{n+1},\nabla v_h\big) +\big(u(\phi_h^n)s_h^{n+1},v_h\big) ,\\
\big(\psi_h^{n+1},v_h\big) &= \big( \phi_h^{n+1},v_h\big)-\big(\nabla\phi_h^{n+1},\nabla v_h\big),\\
s_h^{n+1}-s_h^{n} &=
\frac{1}{2}\big(u(\phi_h^n),\phi_h^{n+1}-\phi_h^{n}\big).
\end{aligned}
\end{equation}
The corresponding error analysis is given by the Theorem \ref{thm3}.
\begin{theorem}\label{thm3}
Let $\phi$ and $\phi_h^n$ be solutions of Eqns.\,\eqref{eq:LB-CH}
and Eqns.\,\eqref{eq:LB-CH-discrete}, respectively. Assume $\phi^0 \in H^2$. In addition, we assume that
$$\phi \in  C(0,T;H^2),\quad \phi_t \in  L^{\infty}(0,T;L^2) \cap L^2(0,T,H^1) \cap C(0,T;H^2),\quad \phi_{tt} \in  L^2(0,T;H^{-1}),$$
and 
$$\psi \in  C(0,T;H^2),\quad  \psi_t \in C(0,T;H^2),\quad \varphi
\in  C(0,T;H^2),\quad  \varphi_t \in C(0,T;H^2).$$
For $0 < n < T/\Delta t$, such that
\begin{equation}\label{eq:eest}
\|e_{\phi,h}^n\|^2+\frac{\xi^2}{2} \| e_{\psi,h}^{n}\|^2 +(e_{s,h}^n)^2 \le C (K_1 \Delta t^2+ K_2 h^4) .
\end{equation} 
where 
\begin{align}
K_1 &= \|\phi_t(r)\|_{L^2(0,T;H^1)}^4+\|\phi_{tt}(r)\|_{L^2(0,T;H^{-1})}^2,\nonumber\\
K_2 &= \|\varphi\|_{C(0,T;H^2)}^2+\|\varphi_t\|_{C(0,T;H^2)}^2+\|\phi\|_{C(0,T;H^2)}^2\nonumber\\
&\quad+\|\phi_t\|_{C(0,T;H^2)}^2+\|\psi\|_{C(0,T;H^2)}^2+\|\psi_t\|_{C(0,T;H^2)}^2.\nonumber
\end{align}
\end{theorem}

It should be pointed out that the above analysis can be
extended to higher-order SAV discretized schemes easily.

\subsection{Linear system}

Let $\{\eta_i\}_{i=1}^N$ be piecewise linear functions which
form a basis of $V_h$. Taking $\phi_h=\sum\limits_{i=1}^N \phi_i
\eta_i$, $\psi_h = \sum\limits_{i=1}^N \psi_i \eta_i$ and
$v_h=\eta_j$, we obtain the matrix form of Eqns.\,\eqref{eq:LB-AC-RFOSAV-LFEM},  
\begin{subequations}
\begin{equation}
\mathbf{M}  \Phi^{n+1}+\Delta t \xi^2 (\mathbf{M} - \mathbf{A}) \Psi^{n+1}+\Delta t \frac{\mathbf{q}}{2}\mathbf{q}^T \Phi^{n+1} = \mathbf{M}   \Phi^{n}  - \Delta t \mathbf{q}[s^{n}-\frac{1}{2}\mathbf{q}^T \Phi^{n}],
\label{eq:LB-AC-RFOSAV-LFEM-A1}
\end{equation}
\begin{equation}
 \mathbf{M}  \Psi^{n+1} = (\mathbf{M} - \mathbf{A}) \Phi^{n+1},
\label{eq:LB-AC-RFOSAV-LFEM-A2}
\end{equation}
\label{eq:LB-AC-RFOSAV-DFEM}
\end{subequations}
where 
$$\Phi:=\{\phi_1,\cdots,\phi_N\}^T,\Psi:=\{\psi_1,\cdots,\psi_N\}^T,\mathbf{q}:=\{q_1,\cdots,q_N\}^T,$$
and
 $$M_{ij} := \big(\eta_i,\eta_j\big), A_{ij} := \big(\nabla \eta_i,\nabla \eta_j\big), q_j :=\big(u_h^{n},\eta_j\big).$$

\subsection{Calculation procedures}
\label{subsec:CD}

Next we can calculate $\Phi^{n+1}$ once $\Phi^{n}$ is known.
Firstly, we denote
\begin{subequations}
\begin{equation}
\mathbf{c}^n := \mathbf{M}  \Phi^{n}  - \Delta t \mathbf{q}[s^{n}-\frac{1}{2}\mathbf{q}^T \Phi^{n}], 
\end{equation}
\begin{equation}
\mathbf{C} := \mathbf{M}+\Delta t \xi^2 (\mathbf{M}-\mathbf{A})\mathbf{M}^{-1}(\mathbf{M}-\mathbf{A}).
\end{equation}
\label{eq:def}
\end{subequations}
Secondly, putting Eqn.\,\eqref{eq:LB-AC-RFOSAV-LFEM2} into
Eqn.\,\eqref{eq:LB-AC-RFOSAV-LFEM1} and using Eqn.\,\eqref{eq:def}, we can obtain
\begin{equation}
\mathbf{C} \Phi^{n+1}+\Delta t \frac{\mathbf{q}}{2}\mathbf{q}^T \Phi^{n+1} = \mathbf{c}^n,
\label{eq:nn11}
\end{equation}
Multiplying Eqn.\,\eqref{eq:nn11} by $\mathbf{C}^{-1}$, and taking the inner
product with $\mathbf{q}^T$, we have 

\begin{equation}
 \mathbf{q}^T \Phi^{n+1}+\Delta t \frac{\mathbf{q}^T\mathbf{C}^{-1}\mathbf{q}}{2}\mathbf{q}^T \Phi^{n+1} = \mathbf{q}^T\mathbf{C}^{-1}\mathbf{c}^n,
\label{eq:nn13}
\end{equation}
The above equation \eqref{eq:nn13} can be rewritten as
\begin{equation}
 \mathbf{q}^T \Phi^{n+1} = \frac{\mathbf{q}^T\mathbf{C}^{-1}\mathbf{c}^n}{1+\Delta t \frac{\mathbf{q}^T\mathbf{C}^{-1}\mathbf{q}}{2}}.
\label{eq:nn14}
\end{equation}
Lastly, $\Phi^{n+1}$ can be obtained by Eqn.\,\eqref{eq:nn11}.

\section{Numerical results}
\label{sec:NR}
In this section, we will present a lot of numerical tests to 
indicate the effectiveness and practicability of our proposed method.

\subsection{Numerical validation} 

We consider one-dimensional problem on uniform mesh with
size $h$ on domain $[0,L](L=4\pi)$, 
starting with an initial solution $u_0=\exp(x/L)$. 
The parameters of LB model are
taken as $\xi=1.0, \alpha=-1.0, \gamma=0.2$. 
 
\subsubsection{Energy dissipation}
\begin{figure}[h]
\centering \mbox{ 
\subfigure[]{\includegraphics[width=45mm]{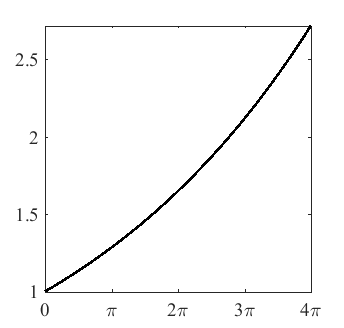}}\quad 
\subfigure[]{\includegraphics[width=45mm]{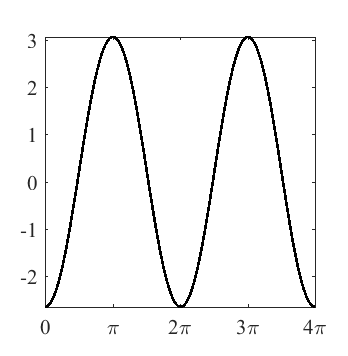}}\quad
\subfigure[]{\includegraphics[width=45mm]{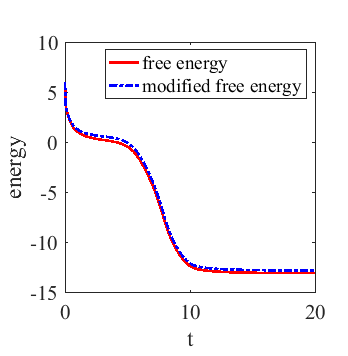}} 
} 
\caption{Lamellar structure: (a) initial phase (when $t=0$); (b) stable phase (when $t=3.22$); (c) free energy}\label{fig:1d-1}
\end{figure}

The first numerical test is done at $h=2^{-8}L, \Delta t=2^{-4}, D_0=16$. 
Fig.\,\ref{fig:1d-1}(a) shows the phase of $\exp(x/L)$.
The stable phase is a lamellar phase including two
periods is represented in Fig.\,\ref{fig:1d-1}(b). 
The number of periods is dominated by the length of computational interval
and the initial value. The energy dissipative property of 
the modified free energy and the original free energy are 
maintained as shown in Fig.\,\ref{fig:1d-1}(c), 
we find that they are close to each other when the
time step length $\Delta t$ is chosen appropriately. 

\subsubsection{Accuracy test}

A series of numerical tests are presented by changing the value
of $\Delta t$ to demonstrate the error order. The numerical
solution with $h=2^{-8}L(L=4\pi), \Delta t=2^{-16}, D_0=25$ is used as a reference
solution to compute error. The $L^2$ numerical errors at $t=2^{-6}$
are shown in Tab.\,\ref{tab:TO}. The first-order accuracy in time
direction is observed in our simulations.

\begin{table}[h]
\centering
\caption{Time errors and convergence rates}\label{tab:TO}
\begin{tabular}{|c|c|c|c|c|c|c|}
\hline
$\Delta t$   &  $\|e_{\phi,h}\|$     &    rate &  $\|e_{\psi,h}\|$     &    rate &  $|e_{s,h}|$     &    rate \\
\hline
$2^{-10}$  & $2.61\text{E-}4$  &$--$  & $2.74\text{E-}3$  &$--$  & $3.04\text{E-}6$  &$--$\\
\hline
$2^{-11}$  & $1.29\text{E-}4$ & $1.01$ & $1.35\text{E-}3$ & $1.02$ & $1.63\text{E-}6$ & $0.90$\\
\hline
$2^{-12}$ & $6.29\text{E-}5$ & $1.04$ & $6.54\text{E-}4$ & $1.05$ & $8.59\text{E-}7$ & $0.92$\\
\hline
$2^{-13}$ &  $2.94\text{E-}5$&  $1.10$ &  $3.05\text{E-}4$ &  $1.10$ &  $4.38\text{E-}7$&  $0.97$\\
\hline
\end{tabular}
\end{table}

To validate the space order of the developed scheme, we design a
set of experiments. By fixing the time step $\Delta t=2^{-12}, D_0=25$ and 
taking the numerical solution with $h=2^{-12}L(L=4\pi)$ as a reference solution, 
we compute every $L^2$ error at $t=2^{-6}$ for different mesh size. 
Tab.\,\ref{tab:SO} illustrates that the second-order accuracy in space 
for the linear FEM.
\begin{table}[h]
\centering
\caption{Space errors and convergence rates}\label{tab:SO}
\begin{tabular}{|c|c|c|c|c|c|c|}
\hline
$h$   &  $\|e_{\phi,h}\|$     &    rate &  $\|e_{\psi,h}\|$     &    rate &  $|e_{s,h}|$     &    rate \\
\hline
$2^{-4}L_0$ &   $1.21\text{E-}2$&  $--$ &  $4.22\text{E-}2$ & $--$ &   $1.07\text{E-}4$ & $--$\\
\hline
$2^{-5}L_0$ &  $3.20\text{E-}3$& $1.92$ &  $1.28\text{E-2}4$ &$1.72$ &  $2.70\text{E-}5$&$1.99$\\   
\hline
$2^{-6}L_0$  &  $8.13\text{E-}4$&  $1.97$&  $3.40\text{E-}3$ & $1.97$ &  $6.72\text{E-}6$ & $2.01$\\
\hline
$2^{-7}L_0$   & $2.03\text{E-}4$ & $2.00$&  $8.62\text{E-}4$ & $1.98$ &   $1.59\text{E-}6$  & $2.08$\\
\hline
\end{tabular}
\end{table}

\subsection{Ordered structures on two-dimensional convex regions} 
\begin{figure}[h]
\centering \mbox{ 
\subfigure[]{\includegraphics[width=45mm]{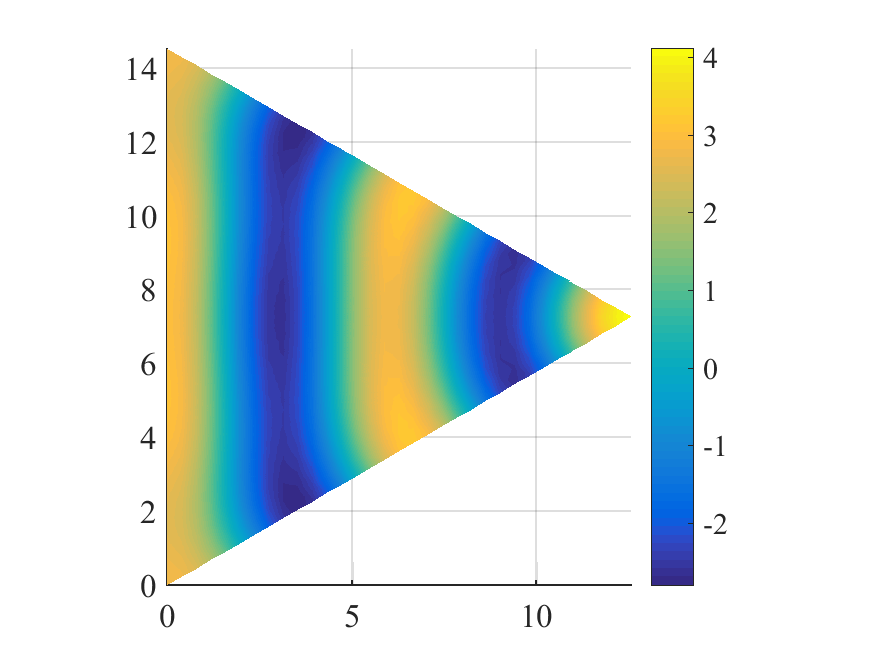}}\quad
\subfigure[]{\includegraphics[width=45mm]{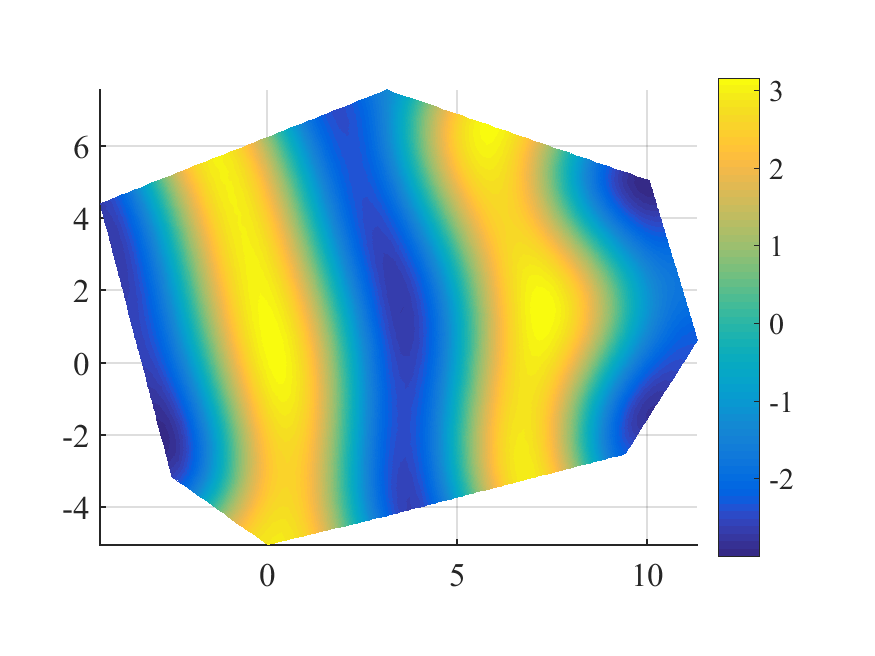}}\quad
\subfigure[]{\includegraphics[width=45mm]{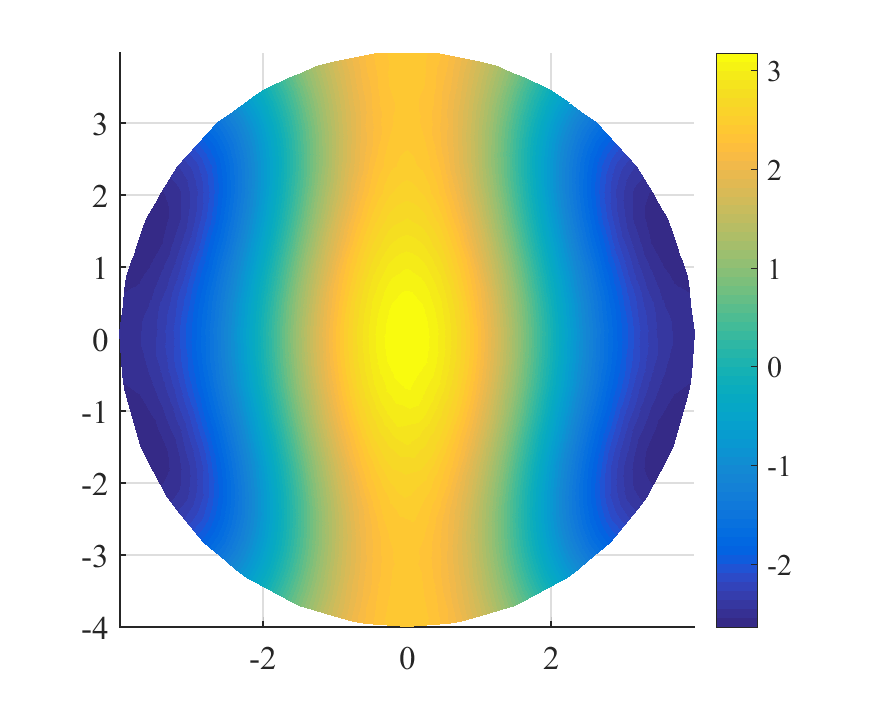}}
}
\centering \mbox{
\subfigure[]{\includegraphics[width=45mm]{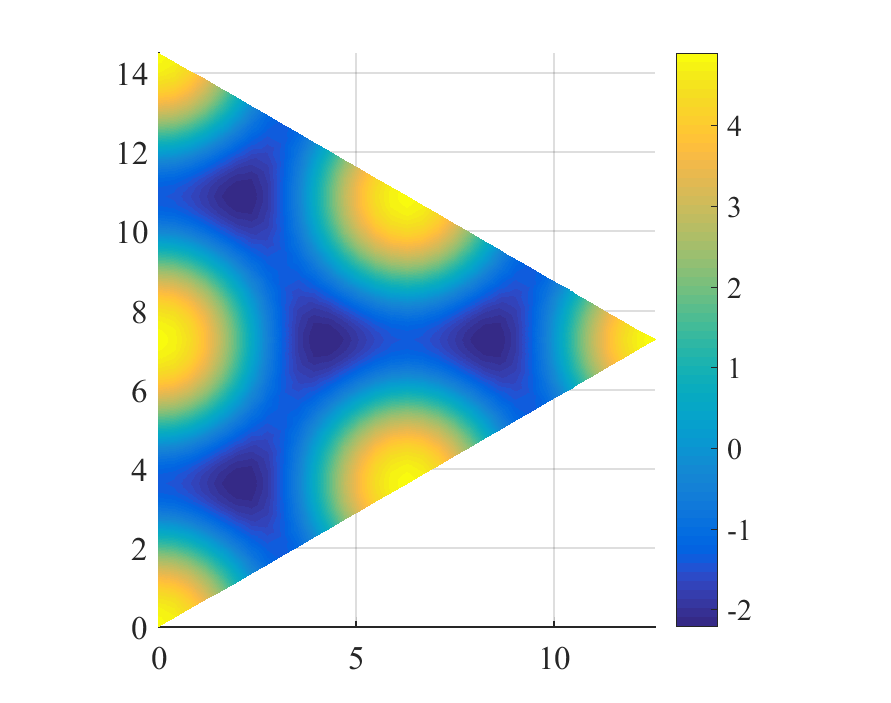}}\quad
\subfigure[]{\includegraphics[width=45mm]{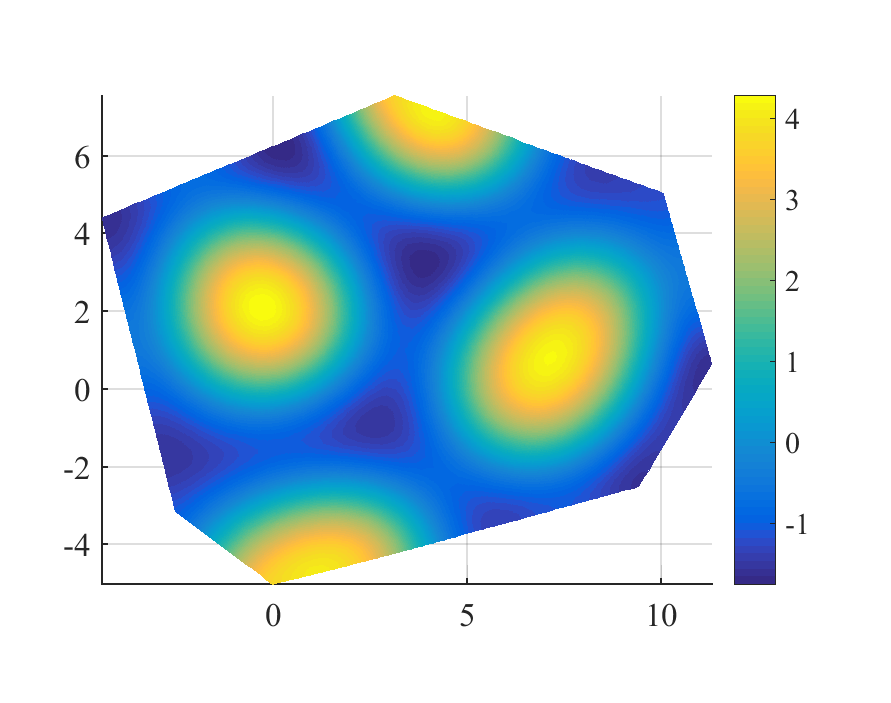}}\quad
\subfigure[]{\includegraphics[width=45mm]{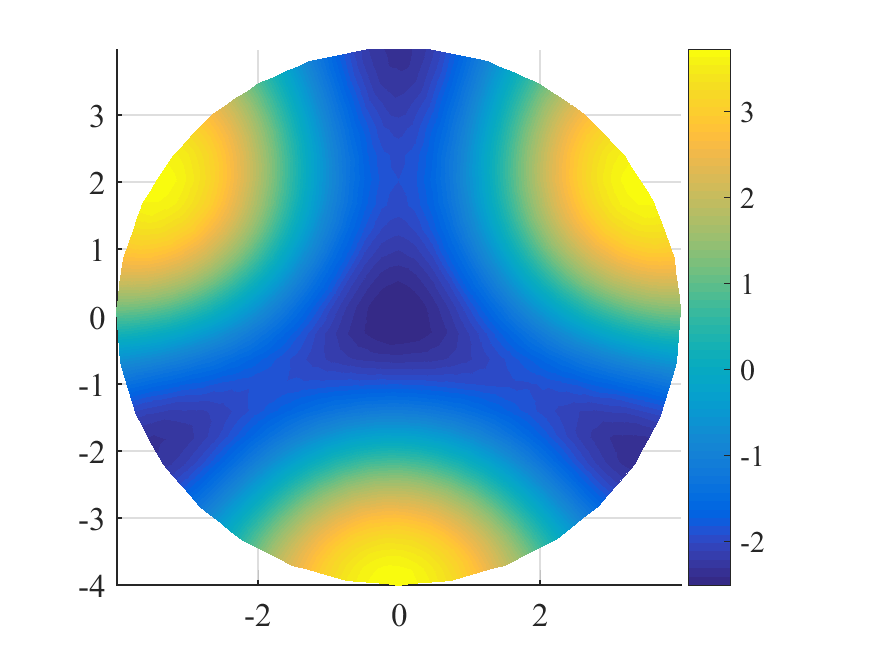}}
} 
\caption{Ordered structures: (a,b,c) lamellar structure; (d,e,f)
hexagonal cylinder structure. (a,d) triangle region; (b,e) hexagon region; (c,f) circle region.}
\label{fig:2d-convex}
\end{figure}

Due to using the splitting method, we only apply our method to convex areas.
We can obtain ordered structures by choosing different parameters
on various two-dimensional convex regions. Some of the results on
different regions, \textit{i.e.} triangle, heptagon, circle regions, are shown in Fig.\,\ref{fig:2d-convex}.
Although we have not shown the diagram of free energy, the energy dissipation
property always holds during time evolution for all tests. These tests are running on the uniform mesh, and the number of corresponding nodes of Fig.\,\ref{fig:2d-convex} (a), (b), (c) are $727$, $1089$, $1503$, respectively.  

\section{Adaptive FEM to PFC model}
\label{sec:AFEM}

In Fig.\,\ref{fig:2d-convex}, we use about $10^3$ mesh nodes to perform numerical simulations on regions containing one or two periodic structures. In fact, we need to do numerical calculations on domains with at least a few dozen periodic structures, such as crystal nucleation and growth. In order to improve computational efficiency, some researchers have worked hard to present various preconditioners\cite{praetorius2011phase,parsons2012numerical,praetorius2015efficient,
farrell2017preconditioner,li2018note,bosch2018preconditioning}. To reduce the computation burden, effective utilization of adaptive mesh refinement to PFC model has been shown\cite{athreya2007adaptive}.  
Here, we shall give a simple adaptive finite element algorithm similar to standard process\cite{verfurth2007adaptive} in Alg.\,\ref{alg:SATA} for PFC model. 

\begin{algorithm}[H]
  \caption{Adaptive finite element algorithm for PFC model}
  \label{alg:SATA} 
\begin{algorithmic}
\REQUIRE  Initial coarse mesh $\mathcal{T}^0$, initial value $\Phi^{0}$, time step $\Delta t$, free energy tolerance $\epsilon_e$, standard deviation tolerance $\epsilon_{\sigma}$.\\
\ENSURE Convergent result $\Phi$ and $\mathcal{E}$  
\STATE{$n:=0$.}
\STATE{Compute free energy $\mathcal{E}^n$.}
\STATE{Set free energy difference $\Delta \mathcal{E} := |\mathcal{E}^{n}|$.}
\WHILE{$\Delta \mathcal{E} > \epsilon_e$}
\STATE{Obtain $\Phi^{n+1}$ on $\mathcal{T}^n$ by calculation procedures \ref{subsec:CD}.} 
\STATE{Calculate $\mathcal{E}^{n+1}$ and $\Delta \mathcal{E} := |\mathcal{E}^{n+1} - \mathcal{E}^{n}|$.}
\STATE{Mesh indicator $\{ \zeta_{\tau}^{n+1} \}_{k\in \mathcal{T}^n}$.}
\STATE{Statistic standard deviation $\sigma$ of $\{\zeta_{\tau}^{n+1}\}$}
\IF{$\sigma > \epsilon_{\sigma}$ }
	\STATE{Coarse mesh.}  
	\STATE{Refine mesh.}
\ENDIF
\STATE{$n:=n+1$.} 
\ENDWHILE 
\end{algorithmic}
\end{algorithm}

\textbf{Remark 5.1}\quad  Since the SAV approach is used in our numerical method, herein the free energy $\mathcal{E}$ in Alg.\,\ref{alg:SATA} is replaced by the modified free energy $\tilde{\mathcal{E}}$.

\textbf{Remark 5.2}\quad 
Our mesh indicator $\zeta_{\tau}$ on each element $k\in \mathcal{T}_h$ have two choices:
\begin{itemize}
\item The standard recovery-type a posteriori error estimator is defined by
\begin{equation}
\zeta_{\tau} := \|\nabla u_h - R_hu_h\|_{0,\tau},
\label{eq:eta}
\end{equation}
\item Gradient estimator is defined as
\begin{equation}
\zeta_{\tau} := \|R_hu_h\|_{0,\tau}.
\label{eq:etatau}
\end{equation}
\end{itemize}

\textbf{Remark 5.3}\quad 
The gradient obtained directly from the numerical solution is discontinuous, the recovered gradient $R_hu_h$ is smoother to be suitable for adaptive indicator. In particular, the SPR technique\cite{zienkiewicz1992superconvergent1,zienkiewicz1992superconvergent2} is used to construct the recovered gradient $R_hu_h$. Certainly, other reconstructed techniques, such as weighted average\cite{bramble1977higher}, PPR\cite{naga2004posterior} or
SCR\cite{huang2010superconvergent}, can be also used to the recovered gradient $R_hu_h$ in our adaptive method.

\textbf{Remark 5.4}\quad 
To balance the distribution of $\zeta_{\tau}$ over the whole region, we introduce standard deviation, which is a concept from statistics. The standard deviation $\sigma$ of $\zeta_{\tau}$ is denoted by
\begin{equation}
\sigma = \sqrt{\sum\limits_{k\in \mathcal{T}}(\zeta_{\tau}-\bar{\zeta})^2/N_{\tau}},
\label{eq:sigma}
\end{equation}
\noindent where $N_{\tau}$ denote the number of elements, $\bar{\zeta}=\sum\limits_{k\in \mathcal{T}} \zeta_{\tau}/N_{\tau}$. 

\textbf{Remark 5.5}\quad 
We use the longest edge bisection algorithm of iFEM\cite{chen2009ifem} to coarse or refine mesh.

\subsection{Efficiency of the gradient estimator for the diblock copolymer system}

In this subsection, we will demonstrate the
gradient estimator is more suitable to the diblock copolymer system
through comparing it with classical $H^1$ error estimator by three examples. 

\subsubsection{Lamellar structure}

\begin{figure}[h]
\centering \mbox{ 
\subfigure[]{\includegraphics[width=45mm]{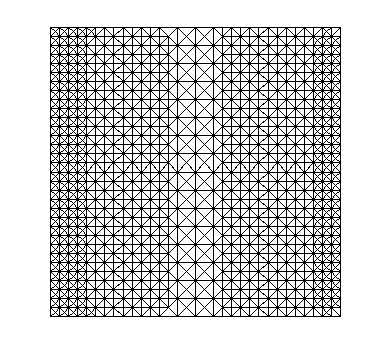}}\quad
\subfigure[]{\includegraphics[width=45mm]{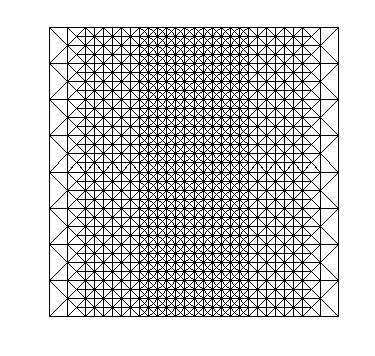}}\quad
\subfigure[]{\includegraphics[width=45mm]{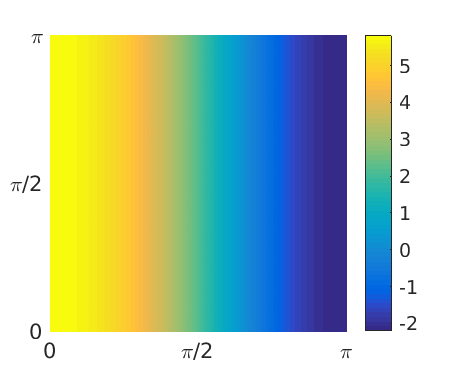}}
}
\caption{Lamellar structure (when $t=4.21$): (a) adaptive mesh by applying $H^1$ error estimator; (b) adaptive mesh by using gradient estimator; (c) convergent phase}
\label{fig:2d-lam}
\end{figure}

We use initial value $u_0(x,y)=\cos(x)$ on the domain
$\Omega=[0,\pi]\times[0,\pi]$. Let $\xi=1.0, \alpha=-1.0, \gamma=0.2,
t=10^{-2},\epsilon_e=10^{-6}, D_0=500$. 
Then we apply the $H^1$ error estimators and 
the gradient estimators in 
our sAFEM to the PFC simulations. 
The lamellar phase can be obtained by these simulations, as shown 
in Fig.\,\ref{fig:2d-lam}(c). 
Fig.\,\ref{fig:2d-lam}(a) depicts the adaptive mesh based on $H^1$ error estimator, which indicates that the mesh near the phase interface has been coarsened. Corresponding, the adaptive mesh by using gradient estimator\ refines near the phase interface (see Fig.\,\ref{fig:2d-lam}(b)). 

\subsubsection{Tetagonal cylinder structure}

\begin{figure}[h]
\centering \mbox{ 
\subfigure[]{\includegraphics[width=45mm]{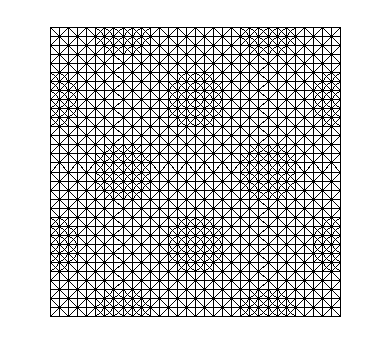}}\quad
\subfigure[]{\includegraphics[width=45mm]{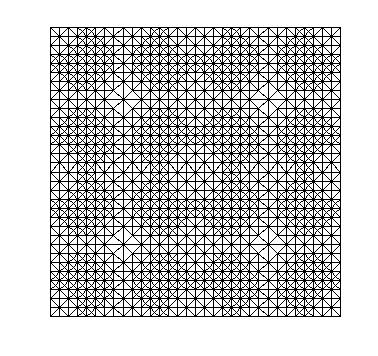}}\quad
\subfigure[]{\includegraphics[width=45mm]{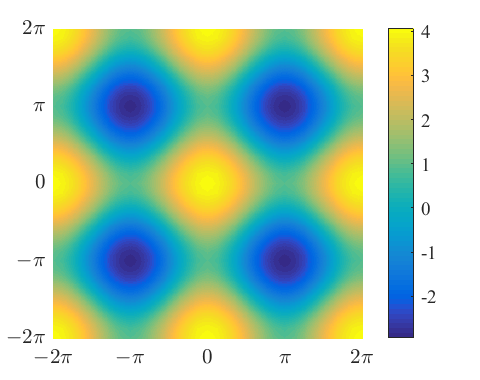}}
}
\caption{Tetagonal cylinder pattern (when $t=12.96$): (a) adaptive mesh by applying $H^1$ error estimator; (b) adaptive mesh by using gradient estimator; (c) convergent morphology}\label{fig:2d-tet}
\end{figure}

We use initial value $u_0(x,y)=\cos(x)+\cos(y)$ on the domain
$\Omega=[-2\pi,2\pi]\times[-2\pi,2\pi]$, then we can obtain the 
tetagonal cylinder structure (see Fig.\,\ref{fig:2d-tet}~(b)) with the
parameters of $\xi=1.0, \alpha=-1.0, \gamma=0.6, \Delta t=10^{-2},  
D_0=500,\epsilon_e=10^{-4}$.
The adaptive mesh, as shown in Fig.\,\ref{fig:2d-tet}~(a) 
comes from sAFEM by using $H^1$ error estimator,
which doesn't guide the adaptive mesh to refine on the phase interface. 
It is also observed that the adaptive meshes by using gradient estimator
are consistent with the density distribution of ordered structures, 
as shown in Fig.\,\ref{fig:2d-tet}~(b) and Fig.\,\ref{fig:2d-tet}\,(c).

\subsubsection{Hexagonal cylinder structure}

\begin{figure}[h]
\centering \mbox{ 
\subfigure[]{\includegraphics[width=45mm]{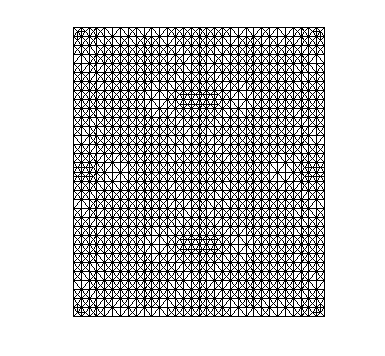}}\quad
\subfigure[]{\includegraphics[width=45mm]{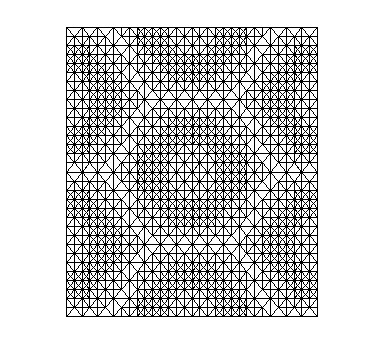}}\quad
\subfigure[]{\includegraphics[width=45mm]{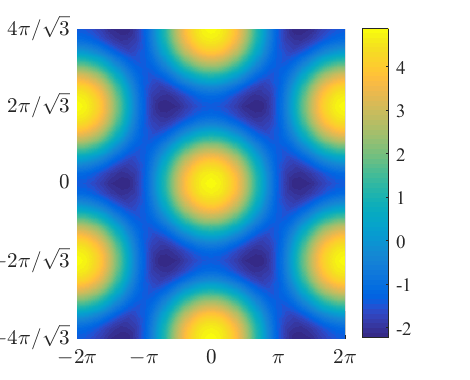}}
}
\caption{Hexagonal cylinder structure (when $t=3.35$): (a) adaptive mesh by applying $H^1$ error estimator; (b) adaptive mesh by using gradient estimator;  (c) free energy}\label{fig:2d-hex}
\end{figure}

When the initial value is chosen as 
$u_0(x,y)=\sum\limits_{j=0}^5 \cos(k_{1, j}x+k_{2, j}y)$,
$k_{1,j}=\cos(j\pi /3)$,
$k_{2,j}=\sin(j\pi /3)$, on the domain of
$\Omega=[-2\pi,2\pi]\times[-4\pi/\sqrt{3},4\pi/\sqrt{3}]$, we can
obtain the hexagonal cylinder structure (see
Fig.\,\ref{fig:2d-hex}~(b)) by sAFEM with $\xi=1.0, \alpha=-1.0,
\gamma=0.8, \Delta t=10^{-2}, D_0=500, \epsilon_e=10^{-4}$. 
The result (Fig.\,\ref{fig:2d-hex}~(b)) based on gradient estimator  
gives a beautiful adaptive mesh, which draws the outline of
hexagonal cylinder structure as shown in Fig.\,\ref{fig:2d-hex}~(c). 
The adaptive meshes by applying $H^1$ error estimator is unsatisfactory.

\textbf{Remark 5.5}\quad 
Zhang \textit{et al.}\cite{zhang2008an} point out that, under the framework of the LB model, when $\xi$ is small, the diblock copolymer system would rather have a macro separation instead of a micro separation. This is not reasonable because only microstructures can form for the diblock copolymer system. Consequently, when simulating the self-assembly behavior of diblock copolymer system, $\xi$ should not be too small. In this paper, the model coefficient $\xi$ is equal to $1$. The ordered structure formed by numerical simulation is smooth in the whole region and the error distribution is nearly uniform. That's why the gradient estimator performs better than the classical $H^1$ error estimator in the above three experiments. Therefore, to better capture the phase interface and its changes of the diblock copolymer system, we use the gradient estimator rather than the classical $H^1$ error estimator.

\subsection{Phase transition} 

In this subsection, we will employ our adaptive method by using gradient estimator 
to simulate the process of phase transition on a square region of 
$[0,6\pi]\times[0,6\pi]$. The initial condition is a mixed state
of hexagonal and lamellar structures. In particular,
$u_0(x,y)$ can be chosen as   
\begin{equation*}
	u_0(x,y) = \left\{
\begin{aligned}
	&6\sin(x+\pi/2),   & x<2\pi,
	\\
	&\sum\limits_{j=0}^5 \cos(k_{1,j}x+k_{2,j}y), 
	k_{1,j}=\cos(j\pi/3),k_{2,j}=\sin(j\pi/3), & x>4\pi,
	\\
	&0, & \mbox{otherwise}.
\end{aligned}
	\right.
\end{equation*}

\begin{figure}[htp]
$\begin{array}{ccc}
$t=0.5$ 
& \tabincell{l}{\includegraphics[width=45mm]{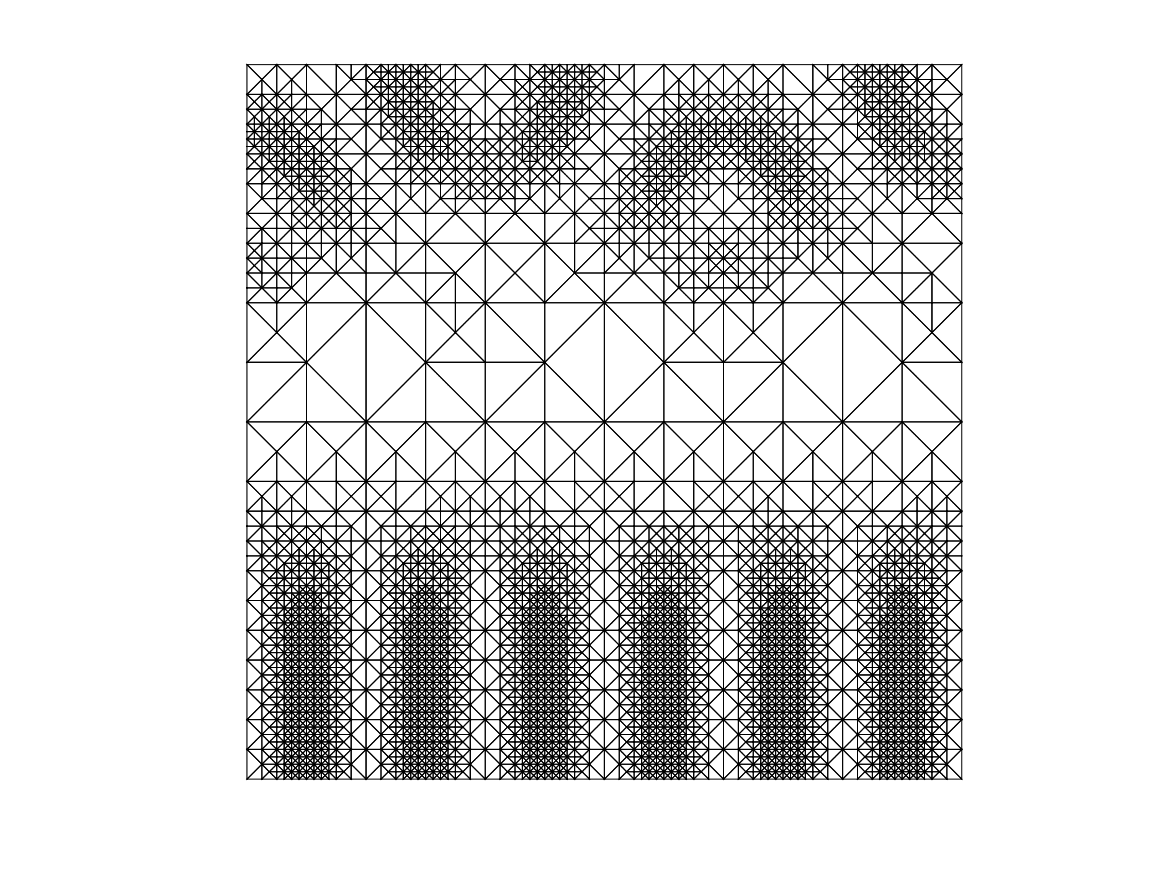}}
& \tabincell{l}{\includegraphics[width=45mm]{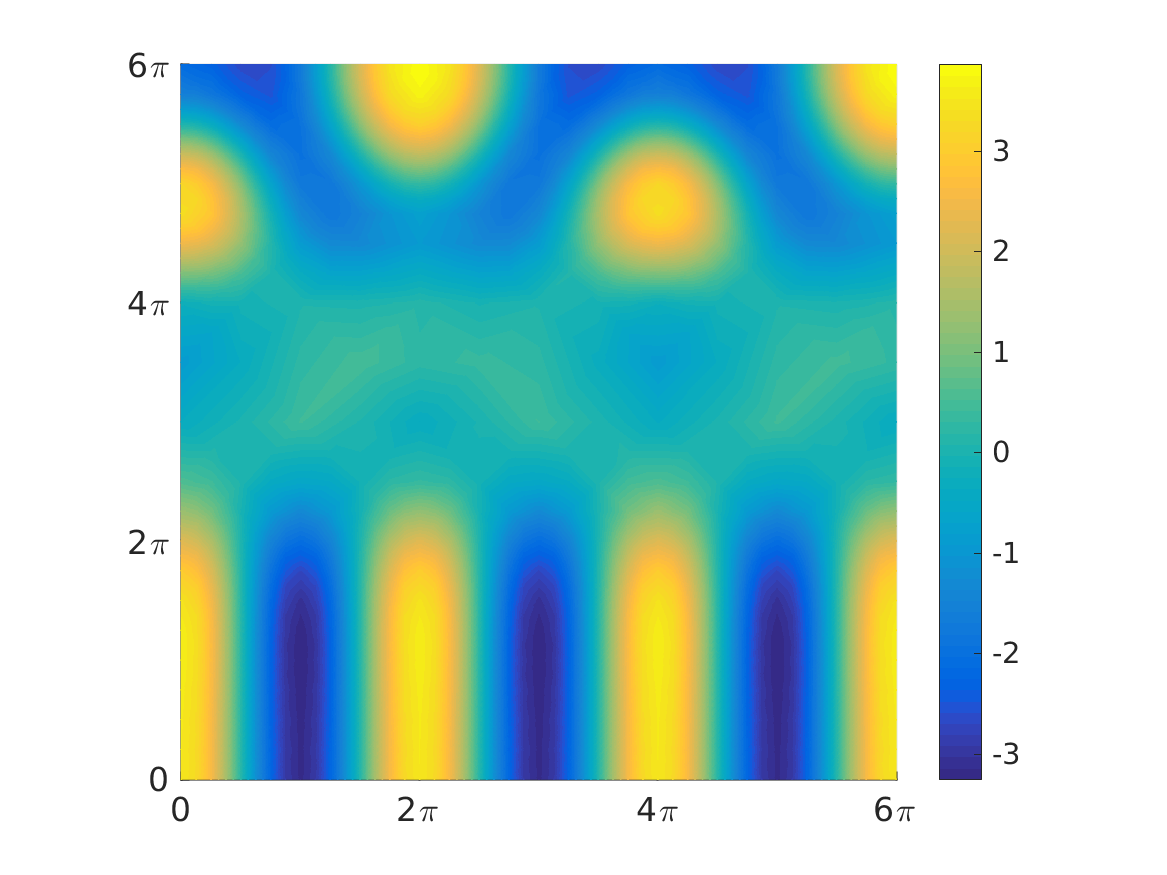}}\\
$t=5$ 
& \tabincell{l}{\includegraphics[width=45mm]{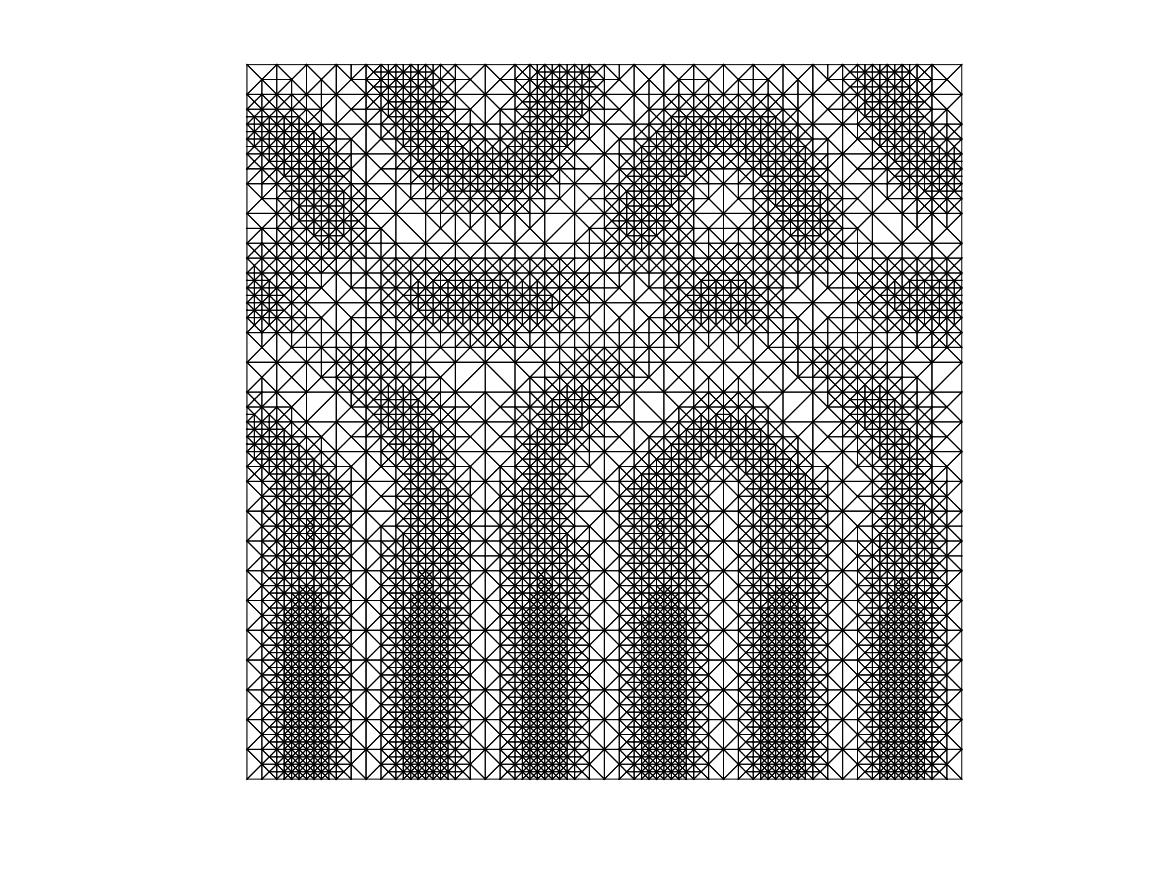}}
& \tabincell{l}{\includegraphics[width=45mm]{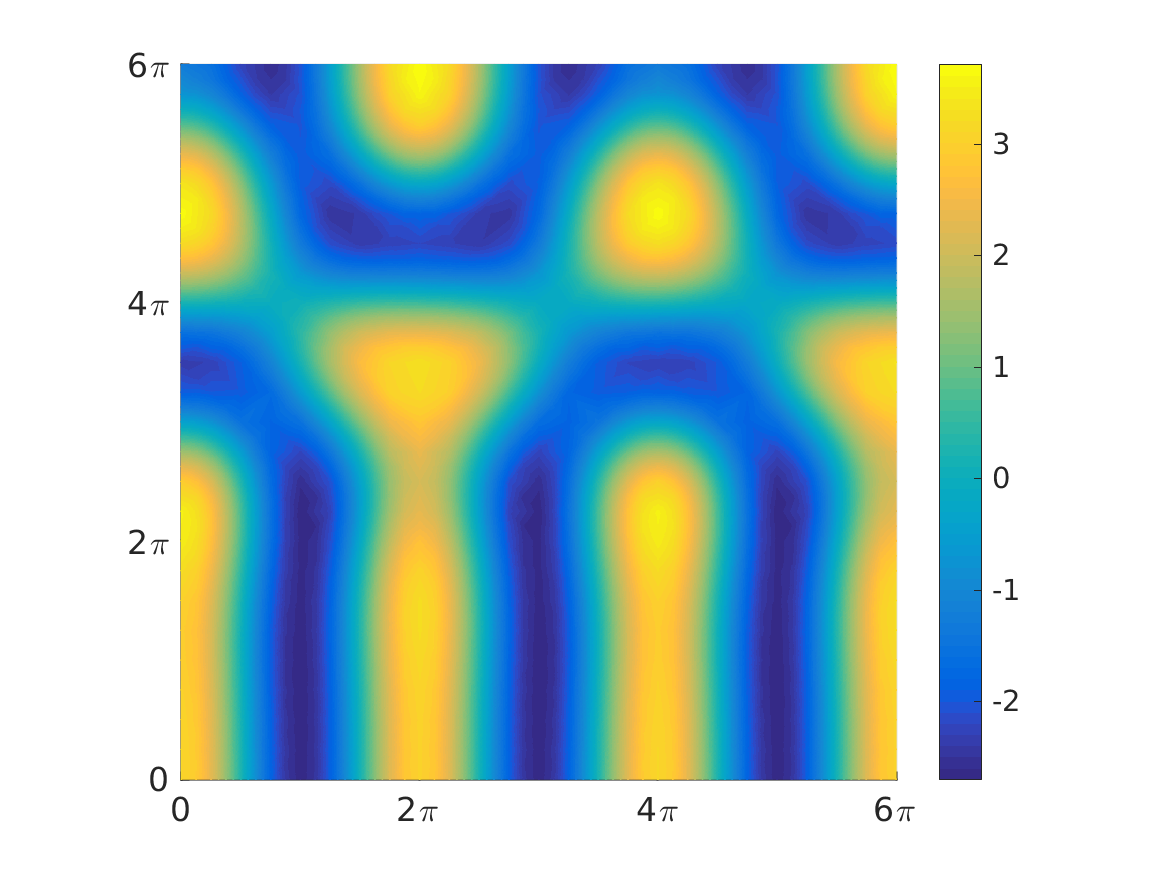}}\\
$t=20$ 
& \tabincell{l}{\includegraphics[width=45mm]{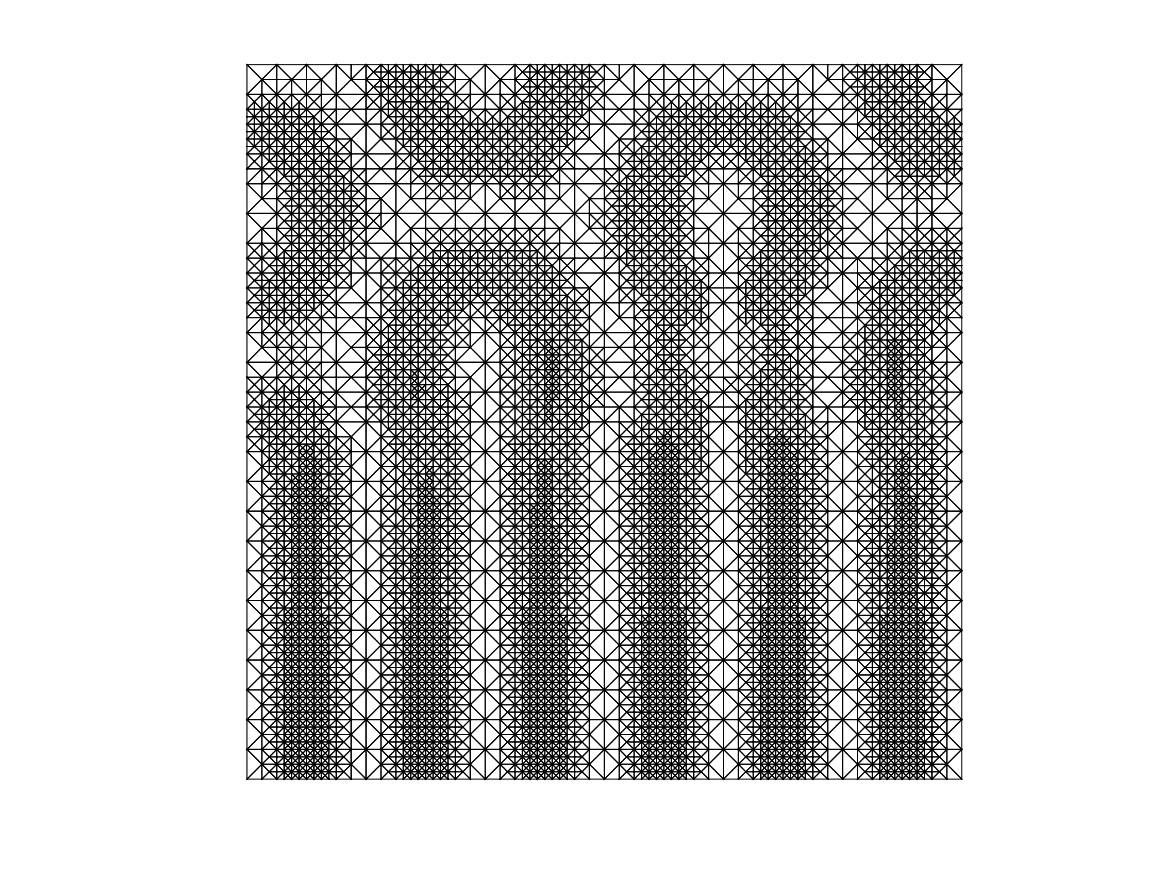}}
& \tabincell{l}{\includegraphics[width=45mm]{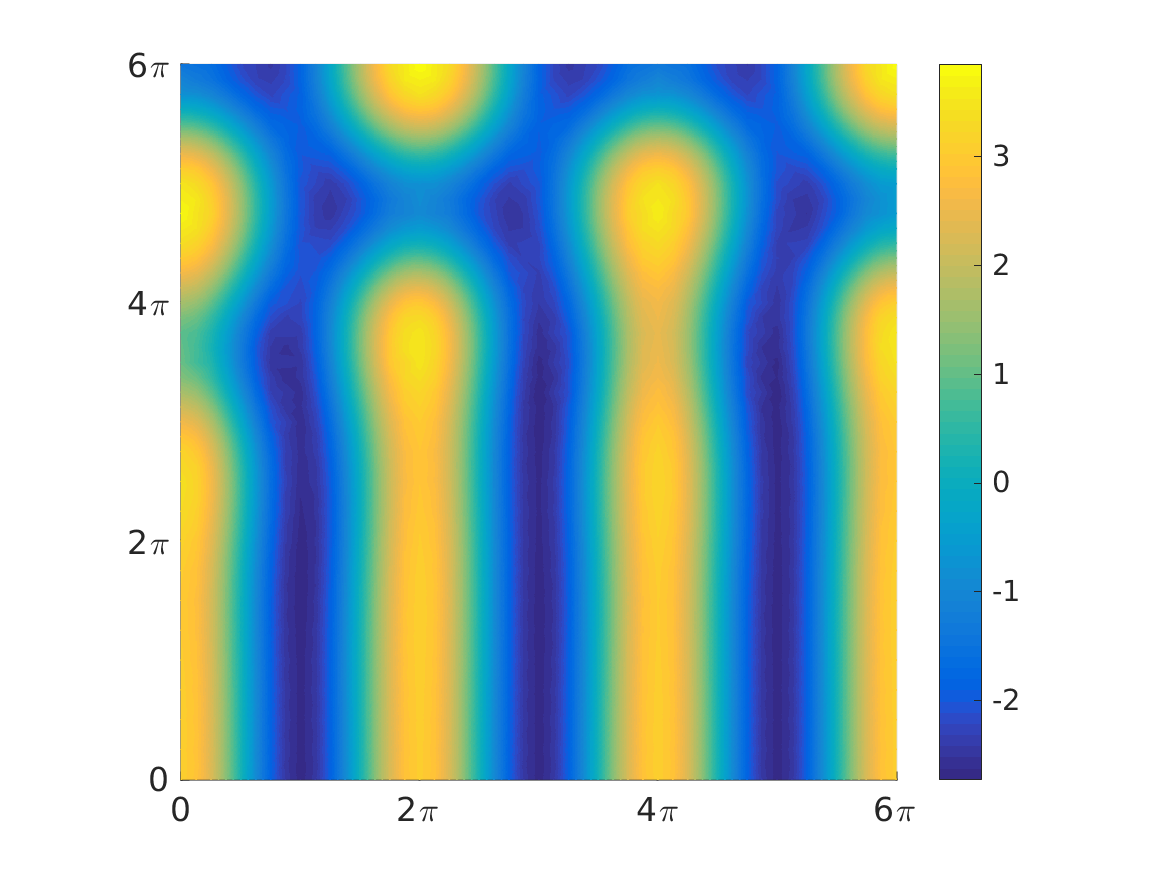}}\\
$t=30$ 
& \tabincell{l}{\includegraphics[width=45mm]{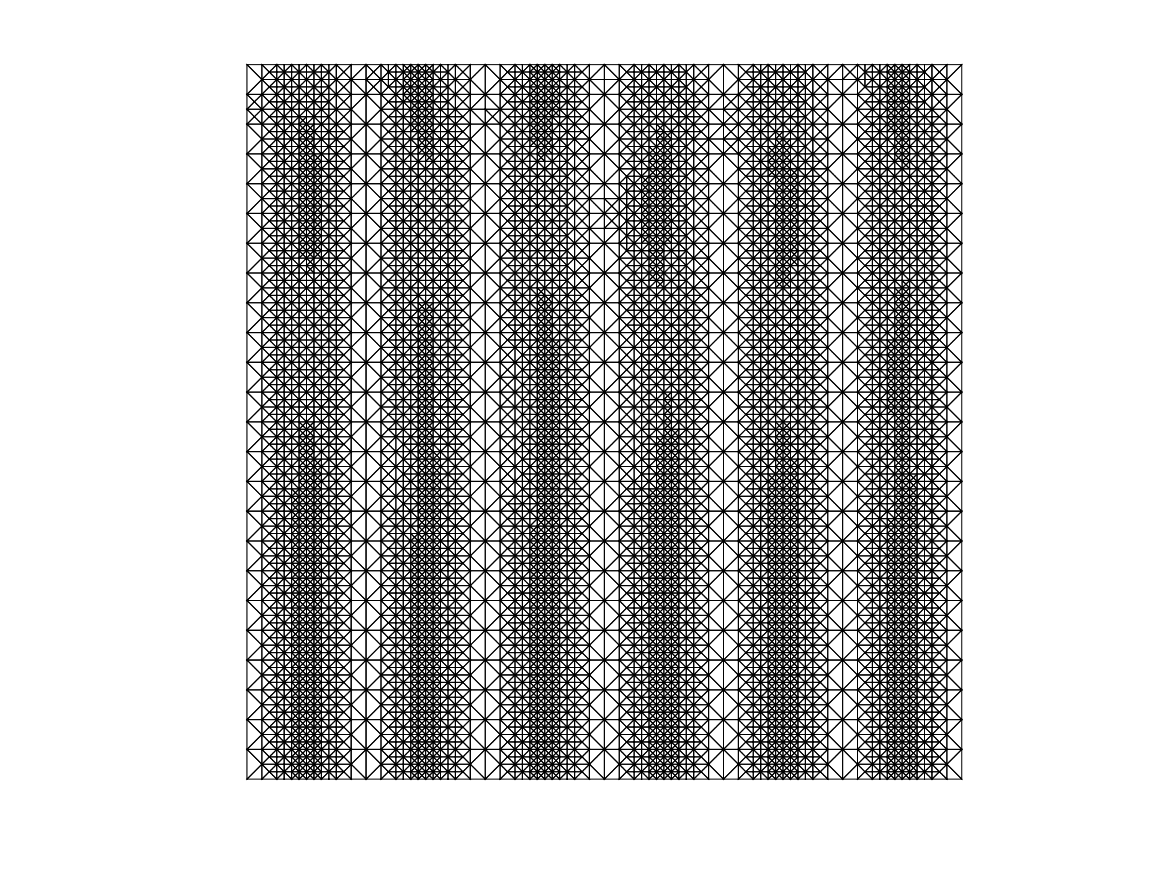}}
& \tabincell{l}{\includegraphics[width=45mm]{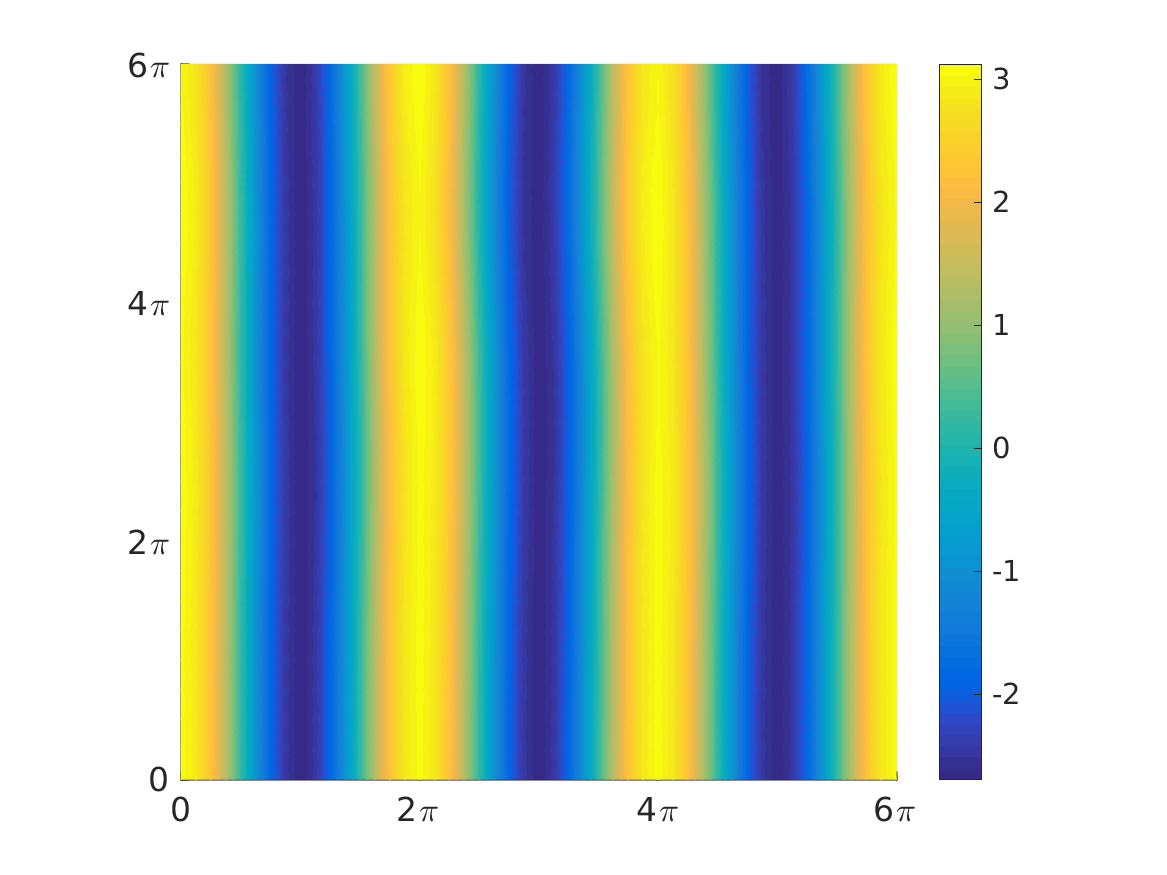}} \\
& (a) & (b) 
 \end{array}$ 
\caption{The dynamical process of phase transition: (a) adaptive meshes; (b) phase evolution} \label{fig:2d-rect-mix} 
\end{figure}

The parameters are $\,\xi=1.0, \alpha=-1.0, \gamma=0.2,
 \Delta t=10^{-2}, D_0=5000,  \epsilon_e=10^{-3},
 \epsilon_{\sigma}=0.05, \theta_r=0.95, \theta_c=0.4$.
Fig.\,\ref{fig:2d-rect-mix} gives the dynamical process.
The left images in Fig.\,\ref{fig:2d-rect-mix} show the adaptive
meshes and the right ones present the evolution process of
morphologies. Due to the lower energy value of the lamellar phase
compared with the hexagonal cylinder pattern, the mixed phase evolves into
lamellar structure as shown in our simulations. More
significantly, our proposed adaptive method can capture 
the interface evolution during the phase transition through the adaptive mesh. 

\section{Conclusions and Outlooks}
\label{sec:CO}
Taking the LB model as an example, we proposed an unconditional
energy stable method, \textit{i.e.} the SAV finite element method, 
to solve the PFC model with Neumann boundary conditions. 
The energy dissipation property of fully discrete scheme has
been proven and error estimate has been derived in theory.
Applying our method, we can effectively simulate the mesoscale
self-assembly in two-dimensional convex geometries. 
To reduce the amount of computing and capture clearly the phase interface, 
a simple adaptive FEM has been present. 
It also should be pointed out that the SAV finite element method and 
adaptive FEM can be improved in several aspects, including:
1) solving the PFC model on concave area;
2) the usage of high-order numerical methods both in time and spatial directions;
3) using parallel technique to solve 3D PFC problems;
4) developing time adaptive methods.



\end{document}